\documentclass{amsart}
\usepackage{epsfig}
\usepackage{graphicx}
\usepackage{amsmath}
\usepackage{amssymb}
\usepackage{amsthm}
\usepackage{graphics}
\usepackage{latexsym}
\begin{document}
\title{Higher Order Modulation Equations for a Boussinesq Equation}
\author{C. Eugene Wayne and J. Douglas Wright}
\date{}

\newtheorem{theorem}{Theorem}
\newtheorem{lemma}{Lemma}
\newtheorem{corollary}{Corollary}
\newtheorem{proposition}{Proposition}
\newtheorem{remark}{Remark}
\newtheorem{hypothesis}{Hypothesis}

\begin{abstract}
In order to investigate corrections to the common KdV approximation to
long waves, we derive modulation  equations for the evolution of long
wavelength initial data for a Boussinesq equation.  The equations
governing the corrections to the KdV approximation are explicitly
solvable and we prove estimates showing that they  do indeed give a
significantly better approximation than the KdV equation alone.  We
also  present the results of numerical experiments which show that the
error estimates we derive are essentially optimal.
\end{abstract}
\maketitle

\section{Introduction}

Modulation, or amplitude, equations are approximate, often explicitly solvable,
model equations derived---usually through asymptotic analysis and the method
of multiple time scales---to model more 
complicated physical situations.  Although these equations have been been used
for over a century, only lately has there been an attempt to rigorously relate
solutions of the modulation equations to the original physical problem.  In particular, 
through the work of Craig \cite{craig:85}, Kano and Nishida \cite{kano.etal:86}, 
Kalyakin \cite{kalyakin:89}, 
Schneider \cite{schneider:98}, Ben Youssef and Colin \cite{youssef.etal:00} 
 and Schneider and Wayne \cite{schneider.etal:00},
\cite{schneider.etal:02}, the validity of Korteweg-de Vries 
(KdV) equations as a leading order approximation
to the evolution of long wavelength water waves and to a number of other dispersive
partial differential equations has been established.  

While the KdV approximation is extremely useful due to its
simplicity and the fact that the KdV equation can be explicitly
solved by the inverse scattering transform, both experimentally and numerically
one  observes departures from the predictions of the KdV equation.  Our goal
in this paper is to derive modulation equations which govern corrections to the KdV
model.  In the present paper we will not work with the full water wave problem
but rather study modulation equations for long wavelength solutions of the
Boussinesq equation:
\begin{equation}
\label{B}
\begin{split}
&\theta_{tt}-\theta_{xx}=(\theta^2)_{xx}+\theta_{ttxx}, \\
&x\in \mathbb{R},\, t\geq 0,\, \theta (x,t)\in \mathbb{R}.
\end{split}
\end{equation}

Our motivation for studying this equation is twofold.  First, the Boussinesq equation
was originally derived as a model equation for water waves, and as our ultimate
goal is to derive corrections to the KdV approximation to water waves, we regard
the study of (\ref{B}) as a useful first step in understanding the much more complicated
water wave situation.  We note that Schneider's analysis of the KdV approximation
for (\ref{B}) in \cite{schneider:98} served as a template for the analysis of the 
water wave problem in \cite{schneider.etal:00}.  

Our second justification for deriving second order modulation  equations for (\ref{B})
is that these modulation equations serve as a sort of normal form for more complicated
PDEs, and as such we expect that the modulation equations which describe corrections
to the KdV approximation to (\ref{B}) will also govern corrections to the KdV approximation
in more complicated situations as well.  Thus the results on existence, uniqueness and other
properties
of the modulation equations we derive in this paper should also be of use in more complicated
situations that we plan to treat in the future.  

We now describe in more detail our results. It is convenient to rewrite
(\ref{B}) as a system of two first order equations.  As in 
\cite{schneider:98} we  introduce new variables,
\begin{equation}\label{uv_def}
\begin{split}
u(x,t)=&\frac{1}{2} (\theta(x,t) - \lambda^{-1} \theta_t(x,t))\\
v(x,t)=&\frac{1}{2} (\theta(x,t) + \lambda^{-1} \theta_t(x,t))
\end{split}
\end{equation}
where $\lambda$ is a skew-symmetric multiplication operator in Fourier space defined 
by $\widehat{\lambda u}=(ik/\sqrt{1+k^2})\hat{u}$.  Note that for $\lambda^{-1} \theta_t$
to be well-defined, we must have $\hat{\theta_t}(0,t)=0$.  That is, the average value
of $\theta_t$ should be zero. We note that for $\theta(x,t)$ a solution of (\ref{B}), we
have that the average value of $\theta_t(x,t)$ is a constant
of the motion. Thus, if the initial condition for $\theta_t$ has
zero average, $\hat{\theta_t}(0,t)=0$ will remain zero for all
time.  Furthermore, as discussed in \cite{schneider:98}, assuming
that the initial condition $\theta_t(x,0)$ has zero average
is not unnatural considering the origin of (\ref{B}).  Thus, we
will make that assumption so that the change of variables
(\ref{uv_def}) is well-defined.

Taking time
derivatives of $u$ and $v$ we find:
\begin{equation}
\label{B2}
\partial _{t}\left( \begin{array}{c}
u\\
v
\end{array}\right) =
\left( \begin{array}{cc}
-\lambda  & 0\\
0 & \lambda 
\end{array}\right) \left( \begin{array}{c}
u\\
v
\end{array}\right) +\frac{1}{2}\left( \begin{array}{c}
-\lambda (u+v)^{2}\\
\lambda (u+v)^{2}
\end{array}\right)
\end{equation}

Not only is (\ref{B2}) convenient from a mathematical point of
view but, as we shall see, $u$ and $v$ have the physical interpretation
of being the left and right moving parts of the solution.

We turn now to the assumptions on the initial conditions
of (\ref{B2}).  The KdV equation is an approximation
of small amplitude and long wavelength motions, 
and thus we will assume that the
initial conditions of (\ref{B}) are of this form.  More precisely, 
fix a constant $C_I > 0$ and assume:
\begin{hypothesis}
\label{hypothesis}
There exist $U_0$, $V_0$ with
$$
\max \left\{ \|U_0\|_{H^s(4) \cap H^{s+9}},\|V_0\|_{H^s(4) \cap H^{s+9}} \right\} < C_I
$$
such that the initial conditions of (\ref{B2}) are of the form:
$$
u(x,0)=\epsilon^2 U_0(\epsilon x) \ ,\ \ 
v(x,0)=\epsilon^2 V_0(\epsilon x)\ ,
$$
where $\epsilon$ is small.
\end{hypothesis}
Here, $H^s(m)=\{f | (1+x^2)^{m/2} f \in H^s\}$.  The norm on this 
weighted Sobolev space
is given by $\|f\|_{H^s(m)} = \|(1+x^2)^{m/2} f \|_{H^s}$.  We use
these spaces because we are interested in solutions which are in some (weak)
sense ``localized''.  In particular, any small perturbation of the known
soliton solutions to the KdV equation will satisfy this localization property.

\begin{remark} We can, of course, recover the initial conditions
for (\ref{B}) from Hypothesis \ref{hypothesis} via (\ref{uv_def}),
and we see that the initial conditions expressed in the $\theta$
variables are also of small amplitude, long-wavelength
form.
\end{remark}

According to the KdV approximation results of \cite{schneider:98},
long wavelength solutions of (\ref{B}) split up into two pieces, one
a right moving wave train and one a left moving wave train.  Each
of these wave trains evolves according to a KdV equation and there is
no interaction between the left and right moving pieces.  One might
expect two types of corrections to such an approximation:
\begin{itemize}
\item corrections due the fact that the left and right moving wave
trains will interact at higher order.
\item corrections due to the fact that even in the case of a purely right (or left)
moving wave train, solutions to the Boussinesq equation are not
exactly described by solutions to the KdV equation.
\end{itemize}
Both of these types of corrections are apparent in our results and in fact
the corrections to the KdV approximation are a sum of solutions of two types
of modulation equations; an inhomogeneous transport equation and a linearized
KdV equation which can be seen (roughly speaking) as reflecting these two
sources of corrections.

To incorporate these two types of corrections, we add to the KdV wavetrains,
which we denote by $U$ and $V$ (since they
represent the leading terms in $u$ and $v$ respectively), 
additional functions $A$ and $B$ and $F$ and $G$.
These functions then satsisfy the modulation equations:
\begin{equation}
\label{KdV}
\begin{split}
\partial_T U =& -\frac{1}{2} \partial_{X_{-}}^3 U - \frac{1}{2} \partial_{X_{-}} U^2\\
\partial_T V =&\phantom{+}  \frac{1}{2} \partial_{X_{+}}^3 V + \frac{1}{2} \partial_{X_{+}} V^2
\end{split}
\end{equation}
\begin{equation}
\label{IH1}
\begin{split}
\partial_\tau A + \partial_X A  =&-\frac{1}{2} 
\partial_{X_+} V^2(X+\tau,\epsilon^2 \tau) - \partial_X U(X-\tau,\epsilon^2 \tau)V(X+\tau,\epsilon^2 \tau)\\
\partial_\tau B - \partial_X B =&\phantom{+}\frac{1}{2} 
\partial_{X_-} U^2(X-\tau,\epsilon^2 \tau) + \partial_X U(X-\tau,\epsilon^2 \tau)V(X+\tau,\epsilon^2 \tau)
\end{split}
\end{equation}
and
\begin{equation}
\label{LIHKdV}
\begin{split}
\partial_T F=&-\partial_{X_{-}} (U F) - \frac{1}{2} \partial^3_{X_{-}} F + J^1\\
\partial_T G=&\phantom{+}\partial_{X_{+}} (V G) + \frac{1}{2} \partial^3_{X_{+}} G + J^2
\end{split}
\end{equation}
where $T=\epsilon^3 t$, $\tau=\epsilon \tau$, $X=\epsilon x$ and $X_{\pm}=X \pm \tau$.

The first of these pairs of equations is simply the KdV approximation.  The second
and third pairs give rise to the corrections to the KdV approximation.  We note
that the terms $J^1$ and $J^2$ which appear in (\ref{LIHKdV}) are inhomogeneous
terms which are made up of a combination of sums and products of the solutions to (\ref{KdV}),
(\ref{IH1}) and their derivatives (see equations (\ref{inhomlkdv1})-(\ref{inhomlkdv2}) below).  

There is some freedom in how we choose the initial data for
the modulation equations.  For simplicity we assume that
$U(X,0)=U_0(X)$, $V(X,0) = V_0(X)$ and
choose zero initial data for (\ref{IH1}) and (\ref{LIHKdV}), i.e.
$A(X,0)=B(X,0)=F(X,0)=G(X,0)=0$.

That the KdV equation has solutions for all times with  initial data
of the type described is well known.  In particular one has (see \cite{schneider.etal:00}):
\begin{theorem}
\label{KdV-solutions} 
Let $s\ge4$.  Then for all $C_0, T_0>0$ there exists $C_1>0$ such that if
$U$, $V$ satisfy (\ref{KdV}) with
initial conditions $U_0$, $V_0$ and 
$$\max \{ \|U_0\|_{H^s(4)\cap H^{s+9}},\|V_0\|_{H^s(4)\cap H^{s+9}}\} < C_0$$
then 
\begin{equation}
\sup_{T\in[0,T_0]} \left\{ \|U(\cdot,T)\|_{H^s(4)\cap H^{s+8}},\|V(\cdot,T)\|_{H^s(4)\cap H^{s+8}} \right\} < C_1
\end{equation}
\end{theorem}
On the other hand it is less clear that solutions of (\ref{IH1}) and (\ref{LIHKdV})
will remain bounded over the very long time scales necessary for the KdV
approximation.  Thus, the first significant technical result of this paper is:
\begin{proposition}
\label{mods_behave}
Fix $T_0 > 0$.  Suppose, $U_0, V_0 \in H^\sigma(4)$ and $U$, $V$, $A$, $B$, $F$ and $G$
satisfy (\ref{KdV})-(\ref{LIHKdV}), then there exists
a constant $C_2$ such that the solutions of (\ref{IH1}) and (\ref{LIHKdV}) satisfy the
following estimates:
\begin{eqnarray*}
\sup_{\tau\in[0,T_0\epsilon^{-2}]} \left\{ \|A(\cdot,\tau)\|_{H^{\sigma-3}},\|B(\cdot,\tau)\|_{H^{\sigma-3}} \right\}
&\le&C_2\\
\sup_{T \in [0,T_0]} \left\{ \|F(\cdot,T)\|_{\tilde{H}}, \|G(\cdot,T)\|_{\tilde{H}} \right\} &\le& C_2
\end{eqnarray*}
where $\tilde{H} = H^{\sigma-5} \cap H^{\sigma-9}(2)$.
\end{proposition}

With this preliminary result in hand, we can now state our principal result.  
\begin{theorem}
\label{main_result}
Fix $T_0,\ C_I > 0,\ \sigma \ge 13$.  
Suppose $U$, $V$, $A$, $B$, $F$ and $G$ satisfy equations (\ref{KdV})-(\ref{LIHKdV}).  
Then there exists $\epsilon_0 > 0$ and $C_F > 0$ such that if
the initial conditions for (\ref{B2}) satisfy Hypothesis \ref{hypothesis}
then for $\epsilon \in (0,\epsilon_0)$, we have that the unique solution to (\ref{B2}) satisfies
$$
\| \bar{u}(\cdot,t) - \bar{w}(\cdot,t) \|_{H^{\sigma-13} \times H^{\sigma-13}} \le C_F \epsilon^{11/2}
$$
for $t \in [0,T_0 \epsilon^{-3}]$,
where,
$$
\bar{w}(x,t)=\epsilon ^{2}\left( \begin{array}{c}
U(X_{-} ,T)\\
V(X_{+} ,T)
\end{array}\right) +\epsilon ^{4}\left( \begin{array}{c}
A(X,\tau )+F(X_{-} ,T)\\
B(X,\tau )+G(X_{+} ,T)
\end{array}\right) 
$$
\end{theorem}
Given this result, and the change of variables
(\ref{uv_def}), we can immediately rewrite this approximation
theorem in terms of the original variables.
Define
\begin{equation*}
\begin{split}
\theta_{app}(x,t)=&\phantom{+}\epsilon^2\left(U(\epsilon(x-t),\epsilon^3 t) + V(\epsilon(x+t),\epsilon^3 t) \right) \\
                  &+\epsilon^4\left(A(\epsilon x, \epsilon t) + B(\epsilon x, \epsilon t) \right)  \\
                  &+\epsilon^4\left(F(\epsilon(x-t),\epsilon^3 t) + G(\epsilon(x+t),\epsilon^3 t) \right). 
\end{split}
\end{equation*}
\begin{corollary}
\label{main_result_theta}
Fix $T_0,\ C_I > 0,\ \sigma \ge 13$.  
Suppose $U$, $V$, $A$, $B$, $F$ and $G$ satisfy equations (\ref{KdV})-(\ref{LIHKdV}).  
Then there exists $\epsilon_0 > 0$ and $C_F > 0$ such that if
the initial conditions for (\ref{B2}) satisfy Hypothesis \ref{hypothesis}
then for $\epsilon \in (0,\epsilon_0)$, the unique solution 
$\theta(x,t)$ to (\ref{B}) satisfies
$$
\| \theta(x,t) - \theta_{app}(x,t)\|_{H^{\sigma-13}} \le C_F \epsilon^{11/2}
$$
for $t \in [0,T_0 \epsilon^{-3}]$.
\end{corollary}

\begin{remark} The initial conditions for equation (\ref{B}) are
obtained from those of (\ref{B2}) simply by inverting the
transformation (\ref{uv_def}).
\end{remark}

The remainder of the paper is devoted to the proof of 
Theorem \ref{main_result}
and Proposition \ref{mods_behave}.  In the next section we give a 
formal derivation
of equations (\ref{KdV})-(\ref{LIHKdV}).  In Section 3 we 
study the existence of 
solutions to equations (\ref{IH1}) and (\ref{LIHKdV}) and prove 
Proposition \ref{mods_behave}.
Section 4 is the technical heart of the paper and contains 
the proof of Theorem \ref{main_result}.
The proof follows the general approach for justifying 
modulation equations laid out
in \cite{kirrmann.etal:92}, but controlling the 
higher order approximation requires
fairly extensive technical modifications.  In Section 5 we 
present the results of a variety of 
numerical computations related to Corollary \ref{main_result_theta}.  
These computations
give insight into several aspects of the second order approximation.  First of
all, it allows us to estimate how large the values of $\epsilon_0$ and $C_F$ in
Corollary \ref{main_result_theta} are.  They also show 
that the order of $\epsilon$
in the error estimates (i.e. $11/2$) is apparently optimal.  Finally, in the concluding
section we discuss other work on second order corrections to the KdV approximation, both rigorous
and non-rigorous, and how it relates to our own results.

\section{Formal Derivation of the Modulation Equations }

One can derive from (\ref{B2}) a system of KdV equations via the method of
multiple time scales  -- this was done in \cite{schneider:98}, for 
example. We extend that calculation in this section
to include the approximating equations for the next order correction.

To derive the modulation equations we first make the {\it Ansatz}  
\begin{equation}
\label{Ansatz}
\left( \begin{array}{c}
u(x,t)\\
v(x,t)
\end{array}\right) =\epsilon ^{2}\left( \begin{array}{c}
U(X_{-} ,T)\\
V(X_{+} ,T)
\end{array}\right) +\epsilon ^{4}\left( \begin{array}{c}
A(X,\tau )+F(X_{-} ,T)\\
B(X,\tau )+G(X_{+} ,T)
\end{array}\right) + O(\epsilon^6)
\end{equation}
where $\tau =\epsilon t$, $T=\epsilon ^{3}t$, $X=\epsilon x$, $X_{-}=X-\tau$, and $X_{+}=X+\tau$.
The two new time variables are the ``multiple time scales'' spoken of earlier.  For convenience
we will also denote
$\bar{u}=(u,v)^{t}$, $\bar{U}=(U,V)^{t}$, $\bar{A}=(A,B)^{t}$ and $\bar{F}=(F,G)^{t}$.

It may seem somewhat
odd that the $O(\epsilon^4)$ correction consists 
of a sum of functions, as opposed to a single function.  
The reason for this is that for our first order approximation terms, $U$ and $V$,
we are assuming that $u$ and $v$ exhibit only unidirectional motion
(right and left, respectively).  
$A$ and $B$,
loosely, correct for the effect of the interaction of right and
left moving waves, and they evolve on the fast time scale, $\tau$. 
There are also unidirectional second order effects, which we represent
with $F$ and $G$.  Their functional form is the same as that of the first
order terms

In a moment, we will insert (\ref{Ansatz}) into (\ref{B2}), but first we compute
the effect of the operator $\lambda$ on long wavelength data. 
Define a function $W$ by $w(x)=W(X)$. We wish to compute $\lambda w(x)$,
and see how that relates to $W(X)$. 
The function $W(X)$ represents solutions of long wavelength and hence
in the Fourier domain we expect the frequency content of these waves to
be concentrated near zero.  Thus, we will (formally) approximate
the effect of $\hat{\lambda}(k)$ by the first few terms of its Maclaurin
series.
\begin{align*}
\lambda w(x) &= \mathfrak{F}^{-1} \left\{\hat{\lambda}(k) \hat{w}(k) \right\}(x)\\
             &= \int e^{ikx} \hat{\lambda}(k) \epsilon^{-1} \hat{W}(k/\epsilon) dk, \quad K=k/\epsilon\\
             &= \int e^{iKX} \hat{\lambda}(\epsilon K) \hat{W}(K) dK\\
             &= \mathfrak{F}^{-1} \left\{\hat{\lambda}(\epsilon K) \hat{W}(K) \right\}(X)\\
             &= \mathfrak{F}^{-1} \left\{\left(\epsilon   (iK)   + \frac{1}{2} \epsilon^3 (iK)^3 
                       + \frac{3}{8} \epsilon^5 (iK)^5  
                       + O(\epsilon^7)\right)\hat{W}(K) \right\}(X)\\
             &=        (\epsilon   \partial_X  + \frac{1}{2} \epsilon^3 \partial^3_X 
         + \frac{3}{8} \epsilon^5 \partial^5_X  
         + O(\epsilon^7))W(X)
\end{align*}
It is important to note that this approximation is only formally good to $O(\epsilon^7)$.

Now we insert this approximation for $\lambda$ 
and the {\it Ansatz} into (\ref{B2}).  This is a necessarily
messy procedure.  To reduce the notation, anything 
formally $O(\epsilon^9)$ or higher is (more or less) 
disregarded.  Also an additional term is added 
to the {\it Ansatz} of the form:
\begin{equation}
\epsilon^6 \bar{S} = \epsilon^6 \left( \begin{array}{c}
S^1(X,\tau)\\
S^2(X,\tau)
\end{array}\right)
\end{equation}
While this term will be treated in much the same way as 
the other terms in the {\it Ansatz}, it should 
be noted that this is not truly part of the next 
order correction.  It will, however, be quite
useful when we prove the approximation is a good one.

We must re-express the partial derivatives in (\ref{B2}) 
in terms of the new coordinates.  By the
chain rule, we have
\begin{equation*}
\partial_t = -\epsilon \partial_{X_{-}} + \epsilon \partial_{X_{+}} + \epsilon \partial_{\tau} 
+ \epsilon^3 \partial_T.
\end{equation*}
Spatial derivatives of terms of the form $f(X_-)g(X_+)$ or $f(X)g(X_\pm)$ are denoted
by $\partial_X$, though all other spatial derivatives are denoted with respect to the appropriate
coordinate.

So we get on the left hand side of (\ref{B2})
\begin{equation}
\label{LHS1}
\begin{split}
\partial_t \left( \begin{array}{c}
                 u(x,t)\\
                 v(x,t) 
                  \end{array}\right)
&=\epsilon^3
\left( \begin{array}{c}
     -\partial_{X_{-}} U(X_{-},T)\\
      \partial_{X_{+}} V(X_{+},T)       
       \end{array}\right)\\
&+\epsilon^5
\left( \begin{array}{c}
      \partial_{T} U(X_{-},T) + \partial_{\tau} A(X,\tau) - \partial_{X_{-}} F(X_{-},T)\\
      \partial_{T} V(X_{+},T) + \partial_{\tau} B(X,\tau) + \partial_{X_{+}} G(X_{+},T)
       \end{array}\right)\\
&+\epsilon^7
\left( \begin{array}{c}
      \partial_{T} F(X_{-},T) + \partial_{\tau} S^1(X,\tau)\\
      \partial_{T} G(X_{+},T) + \partial_{\tau} S^2(X,\tau)
       \end{array}\right)
\end{split}
\end{equation}

Now we must compute the right hand side of (\ref{B2}).  A routine
calculations yields:

\begin{equation}
\begin{split}
\label{THETERROR}
RHS&=\epsilon^3
 \left( \begin{array}{c}
        -\partial_{X_{-}} U \\
         \partial_{X_{+}} V 
        \end{array}\right)\\
 &+\epsilon^5 
 \left( \begin{array}{c}
        -\frac{1}{2} \partial^3_{X_{-}} U - \partial_X A - \partial_{X_{-}} F\\
         \frac{1}{2} \partial^3_{X_{+}} V + \partial_X B + \partial_{X_{+}} G
        \end{array}\right) \\
 &+\epsilon^5
 \left( \begin{array}{c}
        -\frac{1}{2} \partial_{X_{-}} U^2 - \frac{1}{2} \partial_{X_{+}} V^2 - \partial_X (U V)\\
         \frac{1}{2} \partial_{X_{-}} U^2 + \frac{1}{2} \partial_{X_{+}} V^2 + \partial_X (U V)
	\end{array}\right) \\
 &+\epsilon^7
 \left( \begin{array}{c}
         -\frac{1}{2} \partial^3_X A - \frac{1}{2} \partial^3_{X_{-}} F - \frac{3}{8} \partial^5_{X_{-}} U \\
          \frac{1}{2} \partial^3_X B + \frac{1}{2} \partial^3_{X_{+}} G + \frac{3}{8} \partial^5_{X_{+}} V 
        \end{array}\right) \\
 &+\epsilon^7
 \left( \begin{array}{c}
	- \partial_{X} (U A) - \partial_{X_{-}} (U F) - \partial_{X} (U B) - \partial_{X} (U G)\\
	+ \partial_{X} (U A) + \partial_{X_{-}} (U F) + \partial_{X} (U B) + \partial_{X} (U G)
	\end{array} \right) \\        
 &+\epsilon^7
 \left( \begin{array}{c}
 	- \partial_{X} (V A) - \partial_{X} (V F) - \partial_{X} (V B) - \partial_{X_{+}} (V G)\\
 	+ \partial_{X} (V A) + \partial_{X} (V F) + \partial_{X} (V B) + \partial_{X_{+}} (V G)
	\end{array} \right) \\
 &+\epsilon^7
 \left( \begin{array}{c}
        - \frac{1}{4} \partial_{X_{-}}^3 U^2 - \frac{1}{4} \partial_{X_{+}}^3 V^2 
        - \frac{1}{2} \partial_{X}^3 (U V) - \partial_X S^1\\
        + \frac{1}{4} \partial_{X_{-}}^3 U^2 + \frac{1}{4} \partial_{X_{+}}^3 V^2 
        + \frac{1}{2} \partial_{X}^3 (U V) + \partial_X S^2
	\end{array}\right) + O(\epsilon^9)
\end{split}
\end{equation}
So we see that we can satisfy (\ref{B2}) formally to $O(\epsilon^5)$ by taking
\begin{equation}
\begin{split}
\partial_T U =& -\frac{1}{2} \partial_{X_{-}}^3 U - \frac{1}{2} \partial_{X_{-}} U^2\\
\partial_T V =&\phantom{+}  \frac{1}{2} \partial_{X_{+}}^3 V + \frac{1}{2} \partial_{X_{+}} V^2
\end{split}
\tag{\ref{KdV}}
\end{equation}
and
\begin{equation}
\begin{split}
\partial_\tau A + \partial_X A  =&-\frac{1}{2} 
\partial_{X_+} V^2(X+\tau,\epsilon^2 \tau) - \partial_X U(X-\tau,\epsilon^2 \tau)V(X+\tau,\epsilon^2 \tau)\\
\partial_\tau B - \partial_X B =&\phantom{+}\frac{1}{2} 
\partial_{X_-} U^2(X-\tau,\epsilon^2 \tau) + \partial_X U(X-\tau,\epsilon^2 \tau)V(X+\tau,\epsilon^2 \tau)
\end{split}
\tag{\ref{IH1}}
\end{equation}

Equations (\ref{KdV}) are a pair of uncoupled Korteweg-De Vries Equations.  
That their solutions provide the first-order approximation to
long wavelength solutions of (\ref{B}) was proven in \cite{schneider:98}.
Solutions to the KdV equations are known to exist and be bounded over a long
time scale (see Theorem \ref{KdV-solutions} above).

Equations (\ref{IH1}) are
a set of inhomogeneous transport equations driven by the solutions to the KdV equations, for which
we can write down an explicit formula for the solutions to these equations.  

\begin{corollary}
\label{IH1-solutions}
The solutions to equations (\ref{IH1}) are given by
\begin{equation}
\begin{split}
A(X,\tau) =&  \frac{1}{4} V^2(X-\tau,0) - \frac{1}{4} V^2(X+\tau, \epsilon^2 \tau) \\
           & + \alpha(X-\tau,\epsilon^2 \tau) + A_1(X,\tau)\\
B(X,\tau) =&  \frac{1}{4} U^2(X+\tau,0) - \frac{1}{4} U^2(X-\tau, \epsilon^2 \tau) \\
           & +  \beta(X+\tau,\epsilon^2 \tau) + B_1(X,\tau)
\end{split}
\end{equation}
where
\begin{equation}
\begin{split}
A_1(X,\tau) &= \frac{\epsilon^2}{4} \int^\tau_0 \partial_T V^2(X-\tau+2s,\epsilon^2 s) ds  \\
B_1(X,\tau) &= \frac{\epsilon^2}{4} \int^\tau_0 \partial_T U^2(X+\tau-2s,\epsilon^2 s) ds  \\
\end{split}
\end{equation}
\begin{equation}
\begin{split}
\alpha(X_-,T) &= -\epsilon^{-2} \int_0^T \partial_{X_-} \left( U(X_-,s) V(X_- + 2\epsilon^{-2}s,s) \right) ds\\
\beta(X_+,T) &= \epsilon^{-2} \int_0^T \partial_{X_+} \left( U(X_+ - 2\epsilon^{-2}s,s) V(X_+,s) \right) ds
\end{split}
\end{equation}

\end{corollary}

\begin{proof}
The proof follows directly from  Lemmas \ref{transport-derivative} and \ref{cross-term}, which 
appear in the next section.  They can also be verified by inserting the expressions for $A$ and $B$
back into (\ref{IH1}).  Furthermore, as we prove below, in spite of the prefactor
of $\epsilon^{-2}$, $\alpha$ and $\beta$ remain $O(1)$ for all $0 \le T \le T_0$, for any $T_0$.
\end{proof}

The terms of $O(\epsilon^7)$ in (\ref{THETERROR}) give rise to two sets of linear
evolution equations, one for $F$ and $G$ and one for $S^1$ and $S^2$.  Both are 
inhomogeneous systems of equations due to the presence of terms involving $U$, $V$,
$A$, $B$ and their derivatives.  We have some freedom in the way we split up the inhomogeneous
terms between these equations and we attempt to group them in such a way that it
is easy to estimate the resulting solutions over the long time scales relevant for the 
approximation problem.  In particular, we will break $A$ and $B$ up as in the explicit
solutions above.  We have:

\begin{equation}
\begin{split}
\partial_T F =&-\partial_{X_{-}} (U F) - \frac{1}{2} \partial^3_{X_{-}} F + J^1\\
\partial_T G =&\phantom{+}\partial_{X_{+}} (V G) + \frac{1}{2} \partial^3_{X_{+}} G + J^2
\end{split}
\tag{\ref{LIHKdV}}
\end{equation}
where the inhomogeneous terms $J^1$ and $J^2$ are given by:
\begin{equation}
\label{inhomlkdv1}
\begin{split}
J^1(X_-,T)=&-\frac{3}{8} \partial_{X_-}^5 U(X_-,T) - \frac{1}{4} \partial_{X_-}^3 U^2(X_-,T) \\
           &+\frac{1}{4} \partial_{X_-} U^3(X_-,T) - \frac{1}{8} \partial_{X_-}^3 V^2(X_-,0) \\
           &-\partial_{X_-} \left(U(X_-,T)(\frac{1}{4} V^2(X_-,0) + \alpha(X_-,T)) \right)\\
           &-\frac{1}{2} \partial_{X_-}^3 \alpha(X_-,T) 
\end{split}
\end{equation}
\begin{equation}
\label{inhomlkdv2}
\begin{split}
J^2(X_+,T)=&\phantom{+}\frac{3}{8} \partial_{X_+}^5 V(X_+,T) + \frac{1}{4} \partial_{X_+}^3 V^2(X_+,T) \\
           &-\frac{1}{4} \partial_{X_+} V^3(X_+,T) + \frac{1}{8} \partial_{X_+}^3 U^2(X_+,0) \\
           &+\partial_{X_+} \left(V(X_+,T)(\frac{1}{4} U^2(X_+,0) + \beta(X_+,T)) \right)\\
           &+\frac{1}{2} \partial_{X_+}^3 \beta(X_+,T) 
\end{split}
\end{equation}

The additional terms $S^1$ and $S^2$ should satisfy
\begin{equation}
\label{IH2}
\begin{split}
\partial_{\tau} S^1 + \partial_X S^1 = &J^1_{ct} + J^1_{d} + J^1_{sp}\\
\partial_{\tau} S^2 - \partial_X S^2 = &J^2_{ct} + J^2_{d} + J^2_{sp}
\end{split}
\end{equation}
where
\begin{eqnarray*}
J^1_{ct} &=& 
-\partial_X U(X_-,\epsilon^2 \tau) \left(G(X_+,\epsilon^2 \tau) + \beta(X_+,\epsilon^2 \tau) \right)\\
&&-\partial_X U(X_-,\epsilon^2 \tau)\left(\frac{1}{4} U^2(X_+,0) - \frac{1}{4}V^2(X_+,\epsilon^2 \tau)\right)\\
&&- \partial_X \left( V(X_+,\epsilon^2 \tau) F(X_-,\epsilon^2 \tau) \right) 
 - \frac{1}{2} \partial_X^3 \left(U(X_-,\epsilon^2 \tau) V(X_+,\epsilon^2 \tau) \right) \\
&&- \partial_X \left(V(X_+,\epsilon^2 \tau) 
    (\frac{1}{4} V^2(X_-,0) + \alpha(X_-,\epsilon^2 \tau) - \frac{1}{4} U^2(X_-,\epsilon^2 \tau))\right)
\end{eqnarray*}
\begin{eqnarray*}
J^1_{d} &=& 
- \partial_X \left( V(X_+,\epsilon^2 \tau) G(X_+,\epsilon^2 \tau) \right)
 - \frac{1}{8} \partial_X^3 V^2(X_+,\epsilon^2 \tau)+\frac{1}{4}\partial_X V^3(X_+,\epsilon^2 \tau)
\end{eqnarray*}
\begin{eqnarray*}
J^1_{sp} &=& 
- \frac{1}{2}\partial_X^3 A_1(X,\tau) 
 - \partial_X \left((U(X_-,\epsilon^2 \tau)+V(X_+,\epsilon^2 \tau) (A_1(X,\tau) + B_1(X,\tau))\right)\\
&&-\partial_X\left(V(X_+,\epsilon^2 \tau)(\frac{1}{4}U^2(X_+,0)+\beta(X_+,\epsilon^2 \tau))\right)
\end{eqnarray*}
\begin{eqnarray*}
J^2_{ct} &=& 
\partial_X V(X_+,\epsilon^2 \tau) \left(F(X_-,\epsilon^2 \tau) + \alpha(X_-,\epsilon^2 \tau) \right)\\
&&\partial_X(X_+,\epsilon^2 \tau)\left(\frac{1}{4} V^2(X_-,0) - \frac{1}{4}U^2(X_-,\epsilon^2 \tau)\right)\\
&&+ \partial_X \left( U(X_-,\epsilon^2 \tau) G(X_+,\epsilon^2 \tau) \right) 
 + \frac{1}{2} \partial_X^3 \left(U(X_-,\epsilon^2 \tau) V(X_+,\epsilon^2 \tau) \right) \\
&&+ \partial_X \left(U(X_-,\epsilon^2 \tau) 
    (\frac{1}{4} U^2(X_+,0) + \beta(X_+,\epsilon^2 \tau) - \frac{1}{4} V^2(X_+,\epsilon^2 \tau))\right)
\end{eqnarray*}
\begin{eqnarray*}
J^2_{d} &=& 
 \partial_X \left( U(X_-,\epsilon^2 \tau) F(X_-,\epsilon^2 \tau) \right)
 + \frac{1}{8} \partial_X^3 U^2(X_-,\epsilon^2 \tau)-\frac{1}{4}\partial_X U^3(X_-,\epsilon^2 \tau)
\end{eqnarray*}
\begin{eqnarray*}
J^2_{sp} &=& 
 \frac{1}{2}\partial_X^3 B_1(X,\tau) 
 + \partial_X \left(U(X_-,\epsilon^2 \tau)+ V(X_+,\epsilon^2 \tau) (A_1(X,\tau) + B_1(X,\tau))\right)\\
&&+\partial_X\left(U(X_-,\epsilon^2 \tau)(\frac{1}{4}V^2(X_-,0)+\alpha(X_-,\epsilon^2 \tau))\right)
\end{eqnarray*}

The equations (\ref{LIHKdV}) are our second set of modulation equations for the terms
of $O(\epsilon^4)$ in our long wavelength approximation.  Since they are linearized,
inhomogeneous KdV equations, linearized about a KdV solution, they are in principle
explicitly solvable \cite{sachs:83}.  However, the form of the solution that results is
quite complicated (see \cite{sachs:84} and \cite{haragus.etal:98}) 
and thus it requires some effort to show that 
these solutions remain uniformly bounded in the norms which we use to bound the errors.
As we noted above, the functions $S^1$ and $S^2$ do not actually form a part of the 
approximation at $O(\epsilon^4)$; however, we will show that they remain bounded over
the time scales of interest as a part of controlling the error in our approximation.

\section{Estimates on the Solutions to the Modulation Equations}

Before showing that the approximation is a good one, we must first show that the solutions
to the modulation equations are tractable in the their own right.  Keeping in mind that
our goal is show the approximation to (\ref{B2}) is good for a long time, we need to
show that solutions to the modulation equations are bounded on the appropriate time 
scale, that is, for $t\sim O(\epsilon^{-3})$. First we remark on Theorem \ref{KdV-solutions},
above.

Notice that since $T=\epsilon^3 t$, this theorem states that we have bounded solutions of
(\ref{KdV}) for $t\in[0,T_0/\epsilon^3]$, as we had hoped.  Moreover, since the solutions
to (\ref{KdV}) appear in the other modulation equations (often as inhomogeneities), that
they are reasonably smooth and of rapid decay is crucial to showing that the other 
modulation equations are solvable over a long time, and of appropriate size.  
In particular, we will henceforth take $U_0,\ V_0 \in H^{\sigma}(4)$, where $\sigma$ will
be suitably large.
We now state and prove a number of lemmas.

The first set of lemmas concerns the solutions to inhomogeneous transport equations with
zero initial conditions.  From
the method of charateristics, we have explicit formulas for solutions.

\begin{lemma}
\label{transport-naive}
Suppose
$$
\partial_\tau u \pm \partial_X u = f(X,\tau), \quad u(X,0)=0,
$$
with $\|f(X,\tau)\|_{H^{s}} \le C$ for $\tau \in [0,T_0 \epsilon^{-2}]$.  
Then $\|u(\cdot,\tau)\|_{H^s} \le C \epsilon^{-2}$ for $\tau \in [0,T_0 \epsilon^{-2}]$.
\end{lemma}

\begin{proof}
We have 
$$
u(X,\tau)=\int_0^\tau f(X \mp \tau \pm s,s)\ ds. 
$$
The integrand is bounded by $C$ by the Sobolev embedding theorem.  A naive estimate
on the integral proves the result.
\end{proof}

\begin{lemma}
\label{transport-derivative}
Suppose
$$
\partial_\tau u \pm \partial_X u = \partial_X f(X \pm \tau,\epsilon^2 \tau), \quad u(X,0)=0.
$$
Then
\begin{equation}
\label{moc}
u(X,\tau) = \pm \frac{1}{2} \left( f(X \pm \tau, \epsilon^2 \tau) - f(X \mp \tau,0)\right)
           \mp \frac{\epsilon^2}{2}\int^\tau_0 \partial_T f(X \mp \tau \pm 2s, \epsilon^2 s)\ ds
\end{equation}
Also, if $\|f(\cdot,T)\|_{H^{s}} \le C$ and $\|\partial_T f(\cdot,T)\|_{H^{s}} \le C$ for
$T \in [0,T_0]$, then $\|u(\cdot,\tau)\|_{H^{s}} \le C$ for $\tau \in [0,T_0 \epsilon^{-2}]$.
\end{lemma}

\begin{proof}
One can check this result explicitly.  The estimate on the norm follows as in Lemma \ref{transport-naive}.
\end{proof}

\begin{lemma}
\label{cross-term}
Suppose
$$
\partial_\tau u \pm \partial_X u = l(X+\tau,\epsilon^2 t) r(X-\tau,\epsilon^2 \tau),\quad u(X,0)=0.
$$
with $\| l(\cdot,T) \|_{H^s(4)} \le C$ and  $\| r(\cdot,T) \|_{H^s(4)} \le C$ for $T \in [0,T_0]$,
then
$$
u(X,\tau)=\upsilon(X \mp \tau,\epsilon^2 \tau)
$$
with $\| \upsilon(\cdot,T) \|_{H^s(2)} \le C$ for $T \in [0,T_0]$ (that is for $\tau \in [0,T_0 \epsilon^{-2}]$). 
\end{lemma}

\begin{proof}  See appendix.
\end{proof}

\begin{remark} If $l$ and $r$ are taken to be in $H^s(2)$, a similar proof shows that $\upsilon$ is
in $H^s$ over the long time scale.  
\end{remark}

\begin{remark}Since the proof of the Lemma does 
not make explicit use of the slow time scale
dependence of the inhomogeneous factors $l$ and $r$, 
the proof is still valid if the right hand is of the form
$l(X+\tau)r(X-\tau,\epsilon^2 \tau)$, 
$l(X+\tau,\epsilon^2 \tau)r(X-\tau)$ or $l(X+\tau)r(X-\tau)$.
\end{remark}

\begin{remark} A general study of the growth of solutions 
of the transport equation and related linear equations that arise
in the justification of modulation equations was recently
completely by D. Lannes \cite{lannes:02}. 
\end{remark}

With these results, we may now prove the estimate for $A$ and $B$ in Proposition \ref{mods_behave}.  That is,
we have the following:
\begin{corollary}
If $U_0, V_0 \in H^{\sigma}(4)$, and $U$ and $V$ satisfy (\ref{KdV}), then
$$
\sup_{\tau\in[0,T_0\epsilon^{-2}]} \left\{ \|A(\cdot,\tau)\|_{H^{\sigma-3}},\|B(\cdot,\tau)\|_{H^{\sigma-3}} \right\}
\le C.
$$
\end{corollary}

\begin{proof}
From Corollary \ref{IH1-solutions} we know the form of $A$ and $B$.  By Lemma \ref{cross-term}, we have
$\alpha$ and $\beta$ uniformly bounded in $H^{\sigma-1}(2)$ over the long time scale.  Also, by Lemma 
\ref{transport-derivative} we see that $A_1$ and $B_1$ are in the same space as $\partial_T U^2$ and 
$\partial_T V^2$.  $U$ and $V$ satisfy the KdV equations (\ref{KdV}), so we lose three space derivatives
for the one time derivative here.  That is, $A_1$ and $B_1$ are uniformly bounded in $H^{\sigma-3}$ for
the long time scale.
\end{proof}

We
will occasionally be using an alternate, but equivalent norm, on $H^s(2)$.  It is:
\begin{equation*}
|f|_{H^s(2)} = \sum^s_{j=0} \| (1+x^2) \partial^j_x f(x) \|_{L^2}
\end{equation*}
The associated inner product is denoted by $\langle \cdot,\cdot \rangle_{H^s(2)}$

\begin{lemma}  
\label{inner prod 0}
For $s > 3/2$, if $u \in H^{s}$ then  
$$
\langle u \partial_x f,f \rangle_{H^s} \le C |u|_{H^s} |f|^2_{H^s}.
$$
\end{lemma}

\begin{proof}
The proof is similar 
to and simpler than that of Lemma \ref{inner prod 1}, which follows.
\end{proof}

\begin{lemma}
\label{inner prod 1}
For $s > 3/2$, if $u \in H^{s}(2)$ then  
$$
\langle u \partial_x f,f \rangle_{H^s(2)} \le C |u|_{H^s(2)} |f|^2_{H^s(2)}.
$$
\end{lemma}

\begin{proof}  See appendix.
\end{proof}

\begin{lemma}
\label{inner prod 2}
For $f \in H^s(2) \cap H^{s+4}$,
\begin{equation*}
(f,\partial^3_x f)_{H^s(2)} \le C(\|f\|^2_{H^s(2)} + \|f\|^2_{H^{s+4}})
\end{equation*}
\end{lemma}

\begin{proof}  See appendix.
\end{proof}

We may now prove the estimates on $F$ and $G$ in Proposition \ref{mods_behave}.
That is, we have the following Lemma:
\begin{lemma}
\label{LIHKdV-solutions}
If $U_0, V_0 \in H^\sigma(4)$ and $U$, $V$, $A$, $B$, $F$ and $G$ satisfy 
equations (\ref{KdV})-(\ref{LIHKdV}), then $F$ and $G$
satisfy the estimates:
$$
\sup_{T \in [0,T_0]} \left\{ \|F(\cdot,T)\|_{\tilde{H}}, \|G(\cdot,T)\|_{\tilde{H}} \right\} \le C
$$
where $\tilde{H} = H^{\sigma-5} \cap H^{\sigma-9}(2)$.
\end{lemma}

\begin{proof}
The proof follows from Lemmas \ref{inner prod 0} - \ref{inner prod 2} and Gronwall's inequality.    
We show the details for $F$.
The case for $G$ is entirely analogous.  We take the definition of the inner product
on $\tilde{H}$ to be 
$(\cdot,\cdot)_{\tilde{H}} = (\cdot,\cdot)_{H^{\sigma-9}(2)} + (\cdot,\cdot)_{H^{\sigma-5}}$.  

The inhomogeneity $J^1$ is in $H^{\sigma-5}(2)$ (the term $\partial^5_{X_-} U$ causes the loss of
derivatives).  So we
take the inner product of (\ref{LIHKdV}) with $F$ and apply Lemmas 
\ref{inner prod 0} - \ref{inner prod 2} and arrive at:
$$
\partial_T \|F\|_{\tilde{H}}^2 \le C(\|F\|_{\tilde{H}} + \|F\|_{\tilde{H}}^2) \le C(1+\|F\|_{\tilde{H}}^2)
$$
An application of Gronwall's inequality yields, for $T \in [0,T_0]$:
$$
\|F\|_{\tilde{H}}^2(T) \le C T e^{C T}
$$
which concludes the proof of the Lemma and also of Proposition \ref{mods_behave}.
\end{proof}

We now turn our eyes to the set of equations (\ref{IH2}).  
As there are many terms driving these equations, many different 
techniques are used to show that the equations do not blow up over the long time scale.  
We are aided in this task by
the above lemmas, though certain terms will need special consideration.

\begin{lemma}
\label{IH2-solutions}
Suppose $U$, $V$, $A$, $B$, $F$, $G$, $S^1$ and $S^2$ satisfy equations (\ref{KdV})-(\ref{LIHKdV})
and (\ref{IH2}), then $S^1$ and $S^1$
satisfy the estimates:
$$
\sup_{\tau \in [0,T_0 \epsilon^{-2} ]} 
\left\{ \|S^1(\cdot,\tau)\|_{H^{\sigma-10}}, \|S^2(\cdot,\tau)\|_{H^{\sigma-10}} \right\} \le C
$$
\end{lemma}

\begin{proof}
We shall treat the equation for $S^1$ here.  The situation 
for $S^2$ is completely analogous.  Since equations (\ref{IH2}) are
linear, we can consider the inhomogeneity term by term.
First we notice that we can apply Lemma \ref{cross-term} to bound the growth
coming from all terms in $J^1_{ct}$, while Lemma \ref{transport-derivative} suffices
to control all terms coming from $J^1_{d}$.
Thus these terms cause 
no growth over the long time scale.  We now take a moment to discuss the smoothness
of these terms.  
The least smooth term in $J^1_{d}$ is
$\partial_X(V(X+\tau,\epsilon^2 \tau) G(X+\tau,\epsilon^2 \tau))$.  When we apply
Lemma \ref{transport-derivative} we need to examine the smoothness of
$\partial_T(V(X+\tau,\epsilon^2 \tau) G(X+\tau,\epsilon^2 \tau))$.  Now, $G$ is
uniformly bounded in $H^{\sigma-5}$, (from Lemma \ref{LIHKdV-solutions}), and since
$G$ satisfies a linearized KdV equation, we have 
$\partial_T(V(X+\tau,\epsilon^2 \tau) G(X+\tau,\epsilon^2 \tau))$ uniformly bounded
in $H^{\sigma-8}$.  
However, the least smooth term in $J^1_{ct}$ is the term 
$\partial_X (U(X-\tau,\epsilon^2 \tau) G(X+\tau,\epsilon^2 \tau))$, which is 
uniformly bounded in $H^{\sigma-10}(2)$.  Thus at best $S^1$ is uniformly bounded in
$H^{\sigma-10}$.

Each term in $J^1_{sp}$ will require some special consideration.  These terms are:
\begin{eqnarray}
\label{inhih21}
&-&\partial_X( V(X+\tau,\epsilon^2 \tau) \beta(X+\tau,\epsilon^2 \tau))\\
\label{inhih22}
&-&\partial_X( V(X+\tau,\epsilon^2 \tau) U^2(X+\tau,0))\\
\label{inhih23}
&C &\partial^3_X A_1(X,\tau)\\
\label{inhih24}
&C &\partial_X \left((U(X-\tau,\epsilon^2 \tau) + V(X+\tau,\epsilon^2 \tau))(A_1(X,\tau) + B_1(X,\tau))\right)
\end{eqnarray}
Terms (\ref{inhih21}) and (\ref{inhih22}) are treated with by slight variations on Lemmas \ref{transport-derivative}
and \ref{cross-term}.
The technique by which (\ref{inhih23}) and (\ref{inhih24}) are dealt with relies primarily on the prefactor
of $\epsilon^2$ which appears in the definition of the functions
$A_1$ and $B_1$.  Unfortunately, each computation is rather messy.

In the case of the first of these, we apply Lemma \ref{transport-derivative} and get
\begin{equation}
\label{integrals}
\begin{split}
S(X,\tau) &= -\frac{1}{2}(V(X+\tau,\epsilon^2 \tau) \beta(X+\tau,\epsilon^2 \tau) 
                         -V(X-\tau,0) \beta(X-\tau,0) )\\ \nonumber
          &+  \frac{\epsilon^2}{2} \int^\tau_0 \partial_T V(X-\tau + 2s, \epsilon^2 s) \beta(X-\tau+2s,\epsilon^2 s) ds\\ \nonumber
          &+  \frac{\epsilon^2}{2} \int^\tau_0 V(X-\tau + 2s, \epsilon^2 s) \partial_T \beta(X-\tau+2s,\epsilon^2 s) ds
\end{split}
\end{equation}
The first three terms are easily bounded by the techniques discussed previously 
(namely we replace $\partial_T V$ with
the right hand side of the KdV equation and use naive bounds).  
However, when we replace $\partial_T \beta$, we lose the 
prefactor of $\epsilon^2$.  That is, from equation (\ref{new coord trans}) in the proof
of Lemma \ref{cross-term}, we have
$$
\partial_T \beta(X_+,T) = \epsilon^{-2} \partial_X (U(X_+ -2T\epsilon^{-2},T) V(X_+,T)).
$$
We make this substitution into the last term of (\ref{integrals}) to get:
$$
            \frac{1}{2} \int^\tau_0 V(X-\tau + 2s, \epsilon^2 s) \partial_X (U(X-\tau,\epsilon^2 s)
              V(X-\tau+2s,\epsilon^2 s)) ds
$$
Notice that in this integral we have only terms that lie in the weighted Sobolev spaces, and we can use 
the same techniques used in the proof of Lemma \ref{cross-term} to control this term.

The term (\ref{inhih23}) is very nearly of the form needed to apply Lemma \ref{transport-derivative}.  
The only difference
is that there is no dependence on $\epsilon^2 \tau$ in one of the terms.  
The ideas are essentially the same here
as in the proof of Lemma \ref{transport-derivative}.
Consider,
$$
\partial_\tau S + \partial_X S = -\partial_X (V(X+\tau,  \epsilon^2 \tau) U^2(X+\tau,0)).
$$
The solution to this equation is given by:
\begin{align*}
S(X,\tau)=&-\frac{1}{2} \left\{U^2(X+\tau,0) V(X+\tau,\epsilon^2 \tau) - U^2(X-\tau,0) V(X-\tau,\epsilon^2 \tau) \right\} \\
          &+\frac{\epsilon^2}{2} \int^\tau_0 U^2(X - \tau + 2s,0) \partial_T V(X-\tau+2s,\epsilon^2 s) ds
\end{align*}
If one replaces $\partial_T V(X-\tau+2s,\epsilon^2 \tau)$ in the integral by the right hand side of the KdV equation, 
and then takes naive norms, we find that this term is also
controllable.

We now turn our attention to the final two terms which involve the functions
$A_1$ and $B_1$.  The calculations here are quite messy, though the ideas are straightforward.  We replace
$\partial_T V$ with the right hand side of the KdV equation and then apply a number of the same 
techniques used in proving  Lemmas \ref{transport-derivative} and \ref{cross-term}.  The
factor of $\epsilon^2$ present in the definitions of $A_1$ and $B_1$ is crucial.  Consider
\begin{align*}
&\partial_\tau S + \partial_X S \\
=& C \partial_X^3 A_1(X,\tau)\\ 
=&  C \epsilon^2 \partial^3_X   \int^\tau_0 \partial_T V^2(X-\tau+2s,\epsilon^2 s) ds \\ 
=&  C \epsilon^2 \partial^2_X 
     \int^\tau_0 \left( \partial_s (\partial_T V^2(X-\tau+2s,\epsilon^2 s)) 
      - \epsilon^2 \partial_T^2 V(X-\tau+2s,\epsilon^2 s) \right )ds\\
=& C \epsilon^2 \partial_X^2 \left\{\partial_T V^2(X+\tau,\epsilon^2\tau) - \partial_T V^2(X-\tau,0)\right\}\\
 &+C \epsilon^4 \partial^2_X \int_0^\tau \partial_T^2 V^2(X-\tau+2s,\epsilon^2 s)  ds
\end{align*}
For ease of notation, we will let 
$P(X,\tau)= C \epsilon^2 \partial^2_X \int_0^\tau \partial_T^2 V^2(X-\tau+2s,\epsilon^2 s)  ds$.  Notice that
by taking naive estimates on this function, we have that 
$\|P\|_{H^s} \le C$ for $\tau \in [0,T_0 \epsilon^{-2}]$.  Thus we apply 
Lemma \ref{transport-naive} to this equation to find that $S$ is bounded on the long time interval.

In order to deal with (\ref{inhih24}), we will rewrite $\partial_T V^2$ and $\partial_T U^2$. That is, 
\begin{align*}
\partial_T V^2 &= 2 V \partial_T V\\
               &= V (\partial^3_X V + \partial_X V^2)\\
               &= \partial_X \left(V \partial_X^2 V - \frac{1}{2} (\partial_X V)^2 + \frac{2}{3} V^3 \right)\\
               &= \partial_X \tilde{V}.
\end{align*}
Where $\tilde{V}=V \partial_X^2 V - 1/2 (\partial_X V)^2 + 2/3 V^3$.  
A similar calculation yields $\partial_T U^2 = \partial_X \tilde{U}$, where
$\tilde{U}=-U \partial_X^2 U + 1/2 (\partial_X U)^2 - 2/3 U^3$. 
Notice that $\tilde{U}, \tilde{V} \in H^{\sigma-2}$ for $T \in [0,T_0]$, since
they lose at most two derivatives in comparison with $U$ and $V$.  
Similarly, we have $\partial_T \tilde{U} \in H^{\sigma-5}$.  So consider the equation,
\begin{eqnarray*}
& &\partial_\tau S + \partial_X S\\
&=& C\partial_X \left(U(X-\tau,\epsilon^2 \tau) +V(X+\tau,\epsilon^2 \tau)\right)\left(A_1(X,\tau) + B_1(X,\tau)\right)\\ 
&=& C \epsilon^2 \partial_X \Bigl[\Bigr.\left(U(X-\tau,\epsilon^2 \tau)+V(X+\tau,\epsilon^2 \tau)\right)\\
& &\Bigl.\times\int^\tau_0 \partial_T \left(V^2(X-\tau+2s,\epsilon^2 s) 
                   + U^2(X+\tau-2s,\epsilon^2 s)\right)ds\Bigl.\Bigr]\\
&=& C \epsilon^2 \Bigl[\Bigr.\partial_X \left(U(X-\tau,\epsilon^2 \tau)+V(X+\tau,\epsilon^2 \tau)\right)\\ 
& & \times \int^\tau_0\partial_X\left(\tilde{V}(X-\tau+2s,\epsilon^2 s) 
                                    + \tilde{U}(X+\tau-2s,\epsilon^2 s)\right)ds \Bigl.\Bigr]\\
&=& C \epsilon^2 \partial_X \Biggl[\Biggr. \left(U(X-\tau,\epsilon^2 \tau)+V(X+\tau,\epsilon^2 \tau)\right)\\
& & \times \biggl(\biggr. \tilde{V}(X+\tau,\epsilon^2 \tau) - \tilde{V}(X-\tau,0) 
   - \epsilon^2 \int_0^\tau \partial_T \tilde{V} (X-\tau+2s,\epsilon^2 s) ds\\
&& + \tilde{U}(X-\tau,\epsilon^2 \tau) - \tilde{U}(X+\tau,0) 
  - \epsilon^2 \int_0^\tau \partial_T \tilde{U} (X+\tau-2s,\epsilon^2 s) ds \biggl.\biggr) \Biggl.\Biggr]
\end{eqnarray*}
Notice that by taking naive estimates, the terms 
$Q^1(X,\tau) = \epsilon^2 \int_0^\tau \partial_T \tilde{V} (X-\tau+2s,\epsilon^2 s)\ ds$ and 
$Q^2(X,\tau) = \epsilon^2 \int_0^\tau \partial_T \tilde{U} (X+\tau-2s,\epsilon^2 s)\ ds$  
are uniformly
bounded in $H^{\sigma-5}$
over the long time scale.  
Thus we apply Lemma \ref{transport-naive} to the above equation and find 
that this term is well-behaved over the long time scale.
\end{proof}

\section{The Validity of the Approximation}
In this section we set  prove that the approximation 
to a true solution of (\ref{B2}) made by the {\it Ansatz}
is in fact a good one by completing the proof of Theorem \ref{main_result}.

\begin{proof} (of Theorem \ref{main_result}.)

To prove this theorem we shall need a number of lemmas.
\begin{lemma}
\label{op-est 1}
If $\Phi \in H^{s+1}$, then for $\epsilon < 1$,
$$
\| \lambda \Phi(\epsilon \cdot) \|_{H^s} \le C \epsilon^{1/2} \|\Phi\|_{H^{s+1}}.
$$
\end{lemma}

\begin{proof}
The proof here is analogous to the proof of the following lemma.
\end{proof}

\begin{lemma}
\label{op-est 2}
Let $T_1(y)=y$, $T_3(y)=y+1/2y^3$, and $T_5(y)=y+1/2y^3+3/8y^5$.  Then for $j=1,3,5$ if $\Phi(X) \in H^{s+j+2}$ we
have, for $\epsilon < 1$,
$$
\|\lambda\Phi(\epsilon \cdot)-T_j(\epsilon\partial_X)\Phi(\epsilon \cdot)\|_{H^s}\le C \epsilon^{j+3/2} \|\Phi\|_{H^{s+j+2}}
$$
\end{lemma}

\begin{proof} See appendix.
\end{proof}

Now suppose that there is a solution to (\ref{B2}) of the form,
\begin{equation}
\label{assumed-error}
\bar{u}(x,t)=\epsilon^2 \bar{\Psi}(x,t) + \epsilon^{11/2} \bar{R}(x,t)
\end{equation}
where 
\begin{equation}
\epsilon^2 \bar{\Psi}(x,t) = \epsilon^2 \bar{U} + \epsilon^4 (\bar{A}+\bar{F}) + \epsilon^6 \bar{S},
\end{equation}
and $\bar{R} = (R^1(x,t),R^2(x,t))^{t}$.
We consider the term $\bar{R}$ to be the error in our approximation.  Substituting (\ref{assumed-error}) into
(\ref{B2}), we find that $\bar{R}$ must satisfy the equation,
\begin{equation}
\label{error-evolution}
\begin{split}
\partial _{t}\left( \begin{array}{c}
R^1\\
R^2
\end{array}\right) &=
\left( \begin{array}{cc}
-\lambda  & 0\\
0 & \lambda 
\end{array}\right) \left( \begin{array}{c}
R^1\\
R^2
\end{array}\right) 
+\epsilon^2\left( \begin{array}{c}
-\lambda (\Psi^1+\Psi^2)(R^1+R^2)\\
\lambda (\Psi^1+\Psi^2)(R^1+R^2)
\end{array}\right) \\
&+\frac{\epsilon^{11/2}}{2}\left( \begin{array}{c}
-\lambda (R^1+R^2)^{2}\\
\lambda (R^1+R^2)^{2}
\end{array}\right) 
+\epsilon^{-{11/2}} \textrm{Res} [\epsilon^2 \bar{\Psi}]
\end{split}
\end{equation}
where
\begin{equation}
\label{residual}
\textrm{Res} [\epsilon^2 \bar{\Psi}] = -\partial _{t}\left( \begin{array}{c}
\epsilon^2 \Psi^1\\
\epsilon^2 \Psi^2
\end{array}\right) +
\left( \begin{array}{cc}
-\lambda  & 0\\
0 & \lambda 
\end{array}\right) \left( \begin{array}{c}
\epsilon^2 \Psi^1\\
\epsilon^2 \Psi^2
\end{array}\right) 
+\frac{1}{2}\left( \begin{array}{c}
-\lambda (\epsilon^2 \Psi^1+\epsilon^2 \Psi^2)^{2}\\
\lambda (\epsilon^2 \Psi^1+\epsilon^2 \Psi^2)^{2}
\end{array}\right) 
\end{equation}

We have selected our modulation equations precisely so that this term is small.  By taking the time derivative of 
$\bar{\Psi}$ and making then making substitutions from the modulation equations, we find that
\begin{equation}
\label{menace}
\begin{split}
\textrm{Res}[\epsilon^2 \bar{\Psi}]\\
=&
\epsilon^2
\left( \begin{array}{c}
(T_5(\epsilon \partial_{X})  - \lambda) U\\
-(T_5(\epsilon \partial_{X})  - \lambda) V
\end{array}\right) \\
+& \epsilon^4
\left( \begin{array}{c}
(T_3(\epsilon \partial_{X})  - \lambda) \left( A + F + \frac{1}{2}(U+V)^2 \right)\\
-(T_3(\epsilon \partial_{X})  - \lambda) \left( B + G + \frac{1}{2}(U+V)^2\right)
\end{array}\right) \\
+& \epsilon^6
\left( \begin{array}{c}
(T_1(\epsilon \partial_{X})  - \lambda) \left((U+V)(A+F+B+G) + S^1 \right)\\
-(T_1(\epsilon \partial_{X})  - \lambda) \left( (U+V)(A+F+B+G) + S^1 \right)
\end{array}\right) \\
+& \epsilon^8
\left( \begin{array}{c}
- \lambda \left(2(U+V)(S^1+S^2) + (A+F+B+G)^2 \right)\\
  \lambda \left(2(U+V)(S^1+S^2) + (A+F+B+G)^2 \right)
\end{array}\right) \\
+& 2\epsilon^{10}
\left( \begin{array}{c}
- \lambda \left((A+F+B+G)(S^1+S^2)\right)\\
  \lambda \left((A+F+B+G)(S^1+S^2)\right)
\end{array}\right) \\
+& \epsilon^{12}
\left( \begin{array}{c}
- \lambda \left((S^1+S^2)^2\right)\\
  \lambda \left((S^1+S^2)^2\right)
\end{array}\right) 
\end{split}
\end{equation}
While the algebra that goes into showing this is lengthy, it should be noted that this step is accomplished by undoing
to algebra that goes into deriving the modulation equations formally.  

Notice that in the above expression, all functions are of long wavelength form.  Thus we can apply Lemmas \ref{op-est 1}
and \ref{op-est 2} to prove the following result.

\begin{lemma} Under the hypotheses of Theorem \ref{main_result}, the residual satisfies the estimate:
\label{residual estimate}
$$
\sup_{t \in [0,T_0 \epsilon^{-3}]}
\| \textrm{Res} [\epsilon^2 \bar{\Psi}] \|_{H^{\sigma-13} \times H^{\sigma-13}} \le C \epsilon^{17/2}
$$
\end{lemma}

Notice that the loss of three more derivatives is caused 
by the application of Lemma \ref{op-est 2} to
the term in the fourth line of equation (\ref{menace}), 
since $S^j,\ j=1,2$ are uniformly bounded
in $H^{\sigma-10}$. 

We also need the following fact,
\begin{lemma}
\label{following fact}
\begin{equation}
\begin{split}
&\left(R^2-R^1,\lambda[(\Psi^1 +\Psi^2)(R^1+R^2)] \right)_{H^s}\\
 \le&- \left( \partial_t (R^1+ R^2),(\Psi^1 + \Psi^2) (R^1 + R^2) \right)_{H^s} 
 + C \epsilon^{3} \|\bar{R}\|_{H^s \times H^s} \nonumber
\end{split}
\end{equation}
\end{lemma}
\begin{proof}
See appendix.
\end{proof}

We wish to keep the norm of $\bar{R}$ from
growing too much over the long time scale.  That is, if
we can show that $\|\bar{R}\|$ is $O(1)$ for $t \in [0,T_0 \epsilon^{-3}]$,
we will have shown that our approximation is good.  

The first term on the right hand side of equation
(\ref{error-evolution}) will not cause any growth in the norm, since $(f,\lambda f)_{H^s} = 0$. 
The third term has the prefactor of $\epsilon^{11/2}$, which will assist in controlling it, and we
know from the above Lemma \ref{residual estimate} that the residual is small.  

If we tried to control solutions of (\ref{error-evolution}) by applying a
Gronwall type estimate to the time derivative of $(f,f)_{H^s}$, the
second term would result in growth of the norm which would destroy
our estimate over the time scale of interest.  To avoid this problem
we introduce a new energy functional which yields a norm equivalent
to the $H^s \times H^s$ norm, but which does not suffer from this sort of 
uncontrolled growth.

Thus we define
\begin{equation}
\label{new energy}
E^2_s(\bar{R})=\frac{1}{2}\left( \|\bar{R}\|^2_{H^s \times H^s} 
                               + \epsilon^2 (R^1+R^2,(\Psi^1+\Psi^2)(R^1+R^2))_{H^s} \right)
\end{equation}

That this norm is equivalent to the standard norm on $H^s \times H^s$ can be seen by applying
the Cauchy-Schwarz inequality to the inner product, provided that we have 
$\epsilon^2  \|\Psi^1+\Psi^2\|_{H^s} < 1$.  Thus we use without further comment
$$
\frac{1}{C} \| \bar{R} \|_{H^s \times H^s} \le E_s(\bar{R}) \le C \| \bar{R} \|_{H^s \times H^s}.
$$

We now state and prove a useful lemma.

\begin{lemma}
\label{commuter est}
Set $s>0$.  Suppose $f(x),g(x) \in H^s$ and $\gamma(X) \in H^{s+1}$ where $X=\epsilon x$.  Then
$$
| (f(x),\gamma(\epsilon x) g(x))_{H^s}-(g(x),\gamma(\epsilon x) f(x))_{H^s} |
\le C \epsilon \|f\|_{H^s} \|g\|_{H^s} \|\gamma\|_{W^{s,\infty}}.
$$
\end{lemma}

\begin{proof}
See appendix.
\end{proof}

We now have all the tools needed to finish the proof of the theorem.
\begin{eqnarray*}
& &\partial_t E^2_{\sigma-13}(\bar{R})\\
&=& \frac{1}{2} \partial_t \|\bar{R}\|^2_{H^{\sigma-13} \times H^{\sigma-13}} 
 + \frac{\epsilon^2}{2} \partial_t \left( R^1+ R^2,(\Psi^1 + \Psi^2) (R^1 + R^2) \right)_{H^{\sigma-13}}\\
&=& (R^1,\partial_t R^1)_{H^{\sigma-13}} + (R^2,\partial_t R^2)_{H^{\sigma-13}}\\
& &+ \frac{\epsilon^2}{2} \left( \partial_t (R^1+ R^2),(\Psi^1 + \Psi^2) (R^1 + R^2) \right)_{H^{\sigma-13}}\\
& &+\frac{\epsilon^2}{2} \left((R^1+ R^2),(\Psi^1 + \Psi^2) \partial_t (R^1 + R^2) \right)_{H^{\sigma-13}}\\
& &+\frac{\epsilon^2}{2} \left((R^1+ R^2),\partial_t(\Psi^1 + \Psi^2) (R^1 + R^2) \right)_{H^{\sigma-13}}\\
&\le& (R^1,\partial_t R^1)_{H^{\sigma-13}} + (R^2,\partial_t R^2)_{H^{\sigma-13}}\\
& &+ \epsilon^2 \left( \partial_t (R^1+ R^2),(\Psi^1 + \Psi^2) (R^1 + R^2) \right)_{H^{\sigma-13}}\\
& &+ C\epsilon^3 \|\bar{R}\|^2_{H^{\sigma-13} \times H^{\sigma-13}}\\
& &+ C\epsilon^3 \|\bar{R}\|_{H^{\sigma-13} \times H^{\sigma-13}} \|\partial_t\bar{R}\|_{H^{\sigma-13} \times H^{\sigma-13}}\\
&=& \epsilon^2 \left(R^2-R^1,\lambda[(\Psi^1 +\Psi^2)(R^1+R^2)] \right)_{H^{\sigma-13}}\\
 &&+\epsilon^2 \left( \partial_t (R^1+ R^2),(\Psi^1 + \Psi^2) (R^1 + R^2) \right)_{H^{\sigma-13}}\\
 &&+\epsilon^{11/2} \left(R^2-R^1,\lambda[(R^1+R^2)^2] \right)_{H^{\sigma-13}}\\
 &&+\epsilon^{-{11/2}}(\bar{R},\textrm{Res}[\bar{\Psi}])_{H^{\sigma-13} \times H^{\sigma-13}}\\
 &&+ C\epsilon^3 \|\bar{R}\|^2_{H^{\sigma-13} \times H^{\sigma-13}}\\
 &&+ C\epsilon^3 \|\bar{R}\|_{H^{\sigma-13} \times H^{\sigma-13}} \|\partial_t\bar{R}\|_{H^{\sigma-13} \times H^{\sigma-13}}\\
&\le&+ C\epsilon^{{11/2}} \|\bar{R}\|^3_{H^{\sigma-13} \times H^{\sigma-13}}\\
 &&+ C\epsilon^{3} \|\bar{R}\|_{H^{\sigma-13} \times H^{\sigma-13}}\\
 &&+ C\epsilon^3 \|\bar{R}\|^2_{H^{\sigma-13} \times H^{\sigma-13}}\\
&\le& C\epsilon^3 E^2_{\sigma-13}(\bar{R}) 
+ C\epsilon^{3} E_{\sigma-13}(\bar{R}) + C\epsilon^{11/2} E^3_{\sigma-13}(\bar{R})
\end{eqnarray*}

We now state another lemma, 
which proves that the approximation is good over the long time interval.

\begin{lemma}
\label{gron}
Given $C>0$, $T_0>0$, there exists $\epsilon_0>0$ such that if $\epsilon \in (0,\epsilon_0)$
and
\begin{equation}
\dot{\eta}(T) \le C (1 + \eta(T) + \epsilon^{5/2} \eta^{3/2}(T)),\quad \eta(0)=0,
\end{equation}
for $T \in [0,T_0]$ then $\eta(T) \le 2CT_0 e^{2CT_0}$ for $T \in [0,T_0]$.
\end{lemma}

\begin{proof} See appendix.
\end{proof}

We apply Lemma \ref{gron} to the equation for the energy of the remainder and find that
\begin{equation}
\label{energy fact}
E_{\sigma-13}(\bar{R}(\cdot,t)) \le C
\end{equation}
for $t \in [0,T_0 \epsilon^{-3}]$.  Note that given this
{\it a priori} estimate proving the existence and uniqueness
of solutions of  
(\ref{error-evolution}) is a standard exercise. 
  
We now use the equivalence of $E_{\sigma-13}$ to the typical norm
on $H^{\sigma-13}$, equation (\ref{energy fact}), 
and Lemma \ref{op-est 1} to find, for $t\in[0,T_0 \epsilon^{-3}]$:
\begin{equation}
\begin{split}
   \| \bar{u}(\cdot,t) - \bar{w}(\cdot,t) \|_{H^{\sigma-13} \times H^{\sigma-13}}
= &\| \bar{u}(\cdot,t) - \epsilon^2 \bar{\Psi}(\cdot,t) 
                    + \epsilon^6 \bar{S} (\epsilon \cdot,\epsilon t)\|_{H^{\sigma-13} \times H^{\sigma-13}}\\
= &\| \epsilon^{11/2} \bar{R}(\cdot,t)
                    + \epsilon^6 \bar{S}(\epsilon \cdot,\epsilon t) \|_{H^{\sigma-13} \times H^{\sigma-13}}\\ 
\le&C \epsilon^{11/2}
\end{split}
\end{equation}

This completes the proof.
\end{proof}

\section{Some Numerics}

In this section we show the results of some numerical 
simulations.  We performed these numerics to gain insight
into the qualitative nature of the higher order approximation,
to estimate the values of the constants $C_F$ and $\epsilon_0$
that appear in Theorem \ref{main_result}, and to validate 
the results of said theorem.

We will choose the initial conditions of the system
so that we can use the known solitary wave 
solutions to the KdV equation.  We shall
solve the Boussinesq equation (\ref{B2}) 
numerically.  Though techniques are
known for finding explicit solutions to 
the linearized KdV equation (see \cite{haragus.etal:98}
and \cite{sachs:83}), the resulting expressions are quite
complicated, and so we also solve (\ref{LIHKdV}) numerically.

One may wonder why we should even bother computing higher 
order modulation equations, if we have to 
solve them numerically.  In our situation, numerically computing 
solutions to the Boussinesq equation is not
particularly more complicated or time intensive 
than finding solutions to the linearized KdV equations.
However, our goal is to apply these same ideas to derive
corrections to the KdV approximation for the water wave problem,
whose numerical solution is a much more difficult task.
We expect the same modulation equations to hold in these 
more general and complicated systems.  Thus
for the water wave problem,
numerically solving the modulation equations should 
result in a great reduction in the complexity
of the numerics.

The solutions of (\ref{B}) and (\ref{LIHKdV}) are numerically computed using
methods which are largely based around the pseudo-spectral
techniques for Matlab used in \cite{sattinger.etal:98}.  
Since our equations are relatively simple, Matlab, though slower
than other languages (C or Fortran, for example), performs 
adequately rapidly.  The techniques used are largely built around
the use of the fast Fourier transform (FFT) to compute the various
operators and derivatives, and an iterative technique to compute
the nonlinear terms in (\ref{B}) and the term $\partial_X(U F)$
in (\ref{LIHKdV}).  It is implicit in the time step.  

As noted previously, we use the known explicit solutions to (\ref{KdV}).
Where possible, we find explicit solutions for the various terms of $A$ and $B$.
The notable exception to this is in the computation of $\alpha$ and $\beta$, 
which we compute via routine trapezoidal rule techniques.  

We first consider the head-on 
collision of two solitary waves.  
Note that the head-on collision will take place in the 
initial variable $\theta$, not in
either the $u$ or $v$ variables.  This is because 
we have (formally) decomposed
the system into left and right moving waves when we rewrite the
system as (\ref{B2}).  We therefore take initial conditions
such that $U$ and $V$ will evolve as the well-known
$\textrm{sech}$-squared solitary wave solutions to (\ref{KdV}).  
That is we take:
\begin{align*}
u(x,0)&=6\epsilon^2 \textrm{sech}^2(\epsilon x-10)\\
v(x,0)&=6\epsilon^2 \textrm{sech}^2(\epsilon x+10)
\end{align*}
as initial conditions. 

\begin{figure}
\begin{center}
\epsfig{file=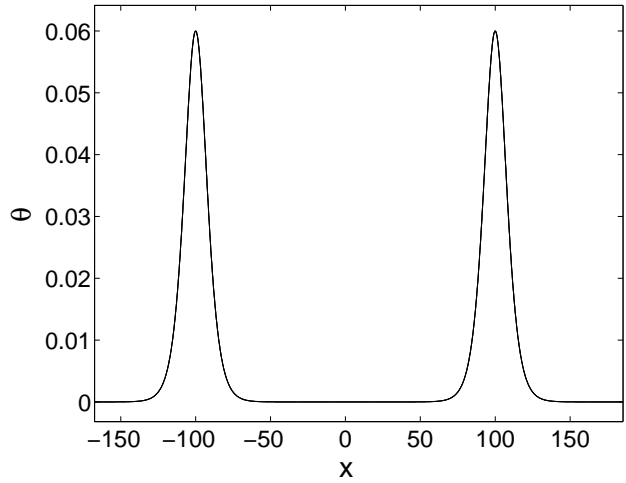, height=2.5in, clip=}
\caption{Initial profile for head-on collision.  $\epsilon = 0.1$}
\label{headoninitial}
\end{center}
\end{figure}

\begin{figure}
\begin{center}
\epsfig{file=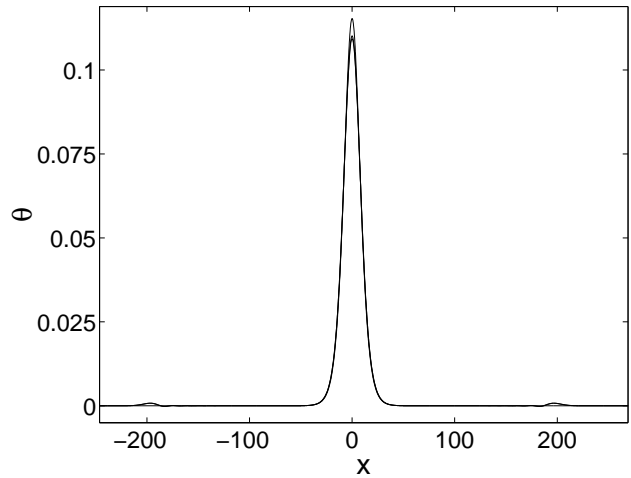, height=2.5in, clip=}
\caption{Head-on collision.  $\epsilon = 0.1$}
\label{headoncollision}
\end{center}
\end{figure}

Figures \ref{headoninitial} and \ref{headoncollision}
show the solution to the Boussinesq equation,
as well as the KdV approximation and the second order correction
on the same plot, at the start and at the collision.  
Here $\epsilon = 0.1$.
We remark on several features of the Boussinesq equation
that are not reflected in the KdV approximation, but are present
in the second order correction.  

\begin{figure}
\begin{center}
\epsfig{file=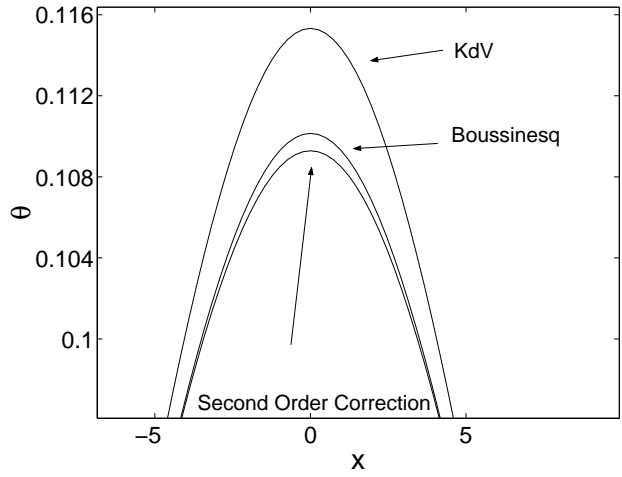, height=2.5in, clip=}
\caption{Close-up of peak of waves, during the head-on collision.  
         $\epsilon = 0.1$}
\label{hocollisionpeak}
\end{center}
\end{figure}

\begin{figure}
\begin{center}
\epsfig{file=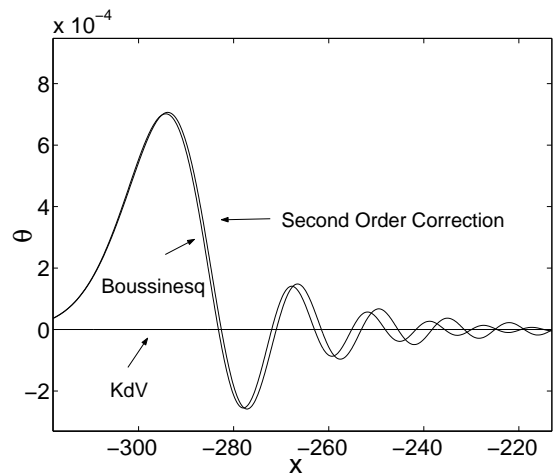, height=2.5in, clip=}
\caption{A ``shadow wave'' and dispersive wave train in the head-on collision.
         $\epsilon = 0.1$}
\label{hocollisiondwave}
\end{center}
\end{figure}

First, in the KdV approximation, during the collision, the 
two waves add in linear superposition (this can be seen as the the KdV
equations evolve independently).  However, the solutions to
the Boussinesq equation do not display this simple linear property
during the collision; the total
height of the wave is slightly less than the sum of the two heights
of the two waves independently.  The second order correction does
a notably better job at displaying this feature 
(see Figure \ref{hocollisionpeak}).
The second feature we notice is the presence of  ``shadow waves'' 
with dispersive wave trains (see Figure \ref{hocollisiondwave}) 
in the solution to (\ref{B}).  These are not present in the 
KdV approximation, but are seen in the second order correction.  

From these pictures, we see that the second order
correction is in fact doing a better job than simply
the KdV approximation alone.  In order to quantify this, 
we computed the solution for variety of 
values of $\epsilon$, and computed the value of the $L^2$ and
$L^\infty$ error of the KdV and second order approximations.  
The time to collision is of $O(\epsilon^{-1})$, and on this time
scale and slightly beyond, the maximum error occurs during the collision.

\begin{figure}
\begin{center}
\epsfig{file=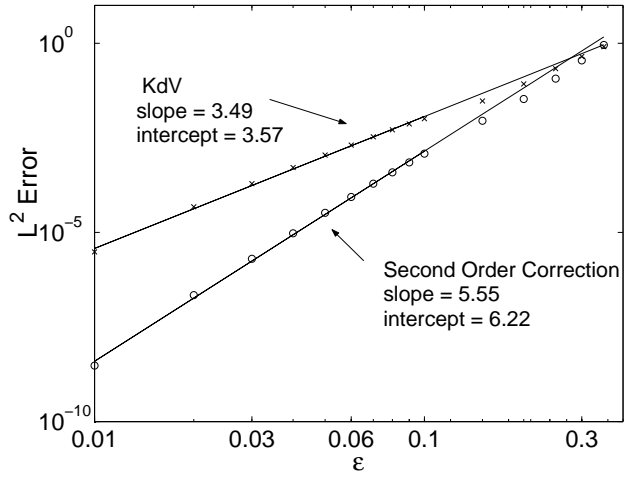, height=2.5in, clip=}
\caption{$\sup \|u - w\|_{L^2}$ vs. $\epsilon$ for head-on collision.}
\label{hol2errvep}
\end{center}
\end{figure}

\begin{figure}
\begin{center}
\epsfig{file=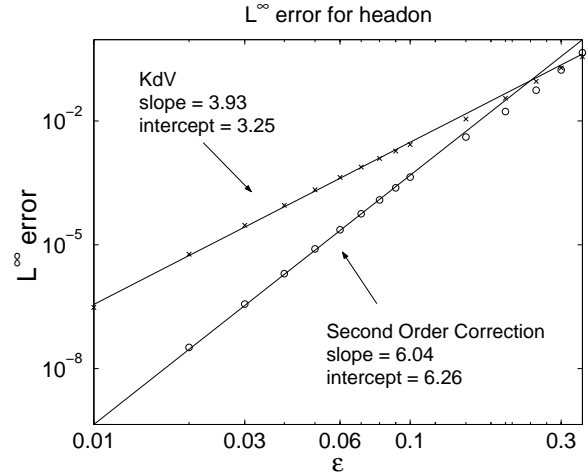, height=2.5in, clip=}
\caption{$\sup \|u - w\|_{L^\infty}$ vs. $\epsilon$ for head-on collision.}
\label{holinferrvep}
\end{center}
\end{figure}

Figures \ref{hol2errvep} and \ref{holinferrvep} 
display log-log plots of the $L^2$ and $L^\infty$ error
versus $\epsilon$ respectively.  The slopes of these lines are 
the order of the correction.  
We note that we have used
only those values of $\epsilon \le 0.1$ in computing these slopes,
as we 
expect the error estimates to hold if $\epsilon$ is sufficiently small.  Table 
\ref{headonslopetable}
summarizes the results.

\begin{table}
\begin{tabular}{|l||c|c|} \hline
                              & $L^2$  & $L^\infty$ \\
\hline\
KdV                           & $3.49$ &  $3.93$\\
\hline
KdV + second order correction & $5.55$ &   $6.04$\\
\hline
\end{tabular}
\caption{Order of the approximation, numerically computed, for the head-on collision}
\label{headonslopetable}
\end{table}

\begin{table}
\begin{tabular}{|l||c|c|} \hline
                              & $L^2$  & $L^\infty$ \\
\hline\
KdV                           & $35.5$ &  $25.8$\\
\hline
KdV + second order correction & $503$ &   $523$\\
\hline
\end{tabular}
\caption{Value of $C_F$, numerically computed, for the head-on collision}
\label{headonconstanttable}
\end{table}

From this we see that our estimate of the error made in
approximating the true solution by the second order
approximation is optimal, in terms of powers of $\epsilon$.
By taking note of the $y$-intercept of these lines, we can
get an estimate on the value of the constant $C_F$ in each case,
see Table (\ref{headonconstanttable}).
Unfortunately, these values of the constant are quite large for the
second order correction.  We also note that it is not so much the
actual value of the leading coefficient that matters, as is the location
(in $\epsilon$) at which the second order correction and the KdV
correction return the same error.  That is, graphically, where the lines
in figures (\ref{hol2errvep}) and (\ref{holinferrvep}) cross.

The next simulation was that of right moving overtaking
waves.  We take initial data such that $U$ will evolve as the famous
two soliton solution to (\ref{KdV}).  Since we are not interested
in left moving waves, we take initial data for $v$ to be zero.  
Note that $v$ does not remain zero, however, due to the coupling.

Unlike 
the previous situation, the time scale of the overtaking wave collision
is $O(\epsilon^{-3})$.  To observe the entirety of the collision,
we take $T_0 = 8$.  We also observe 
that the error in the approximation is largest at the end of the interval
$[0,T_0 \epsilon^{-3}]$.  From the proof of Lemma \ref{gron}, one can see
that as $T_0$ increases, $\epsilon_0$ decreases.  This requires smaller values
of $\epsilon$, which in turn necessitates running the simulation for a longer
period of time.  

\begin{figure}
\begin{center}
\epsfig{file=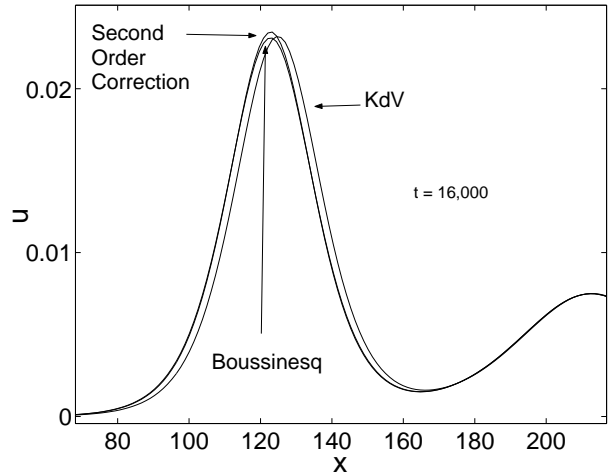, height=2.5in, clip=}
\caption{Overtaking Wave before the Collision, $\epsilon = 0.05$}
\label{overbeforecollision}
\end{center}
\end{figure}

\begin{figure}
\begin{center}
\epsfig{file=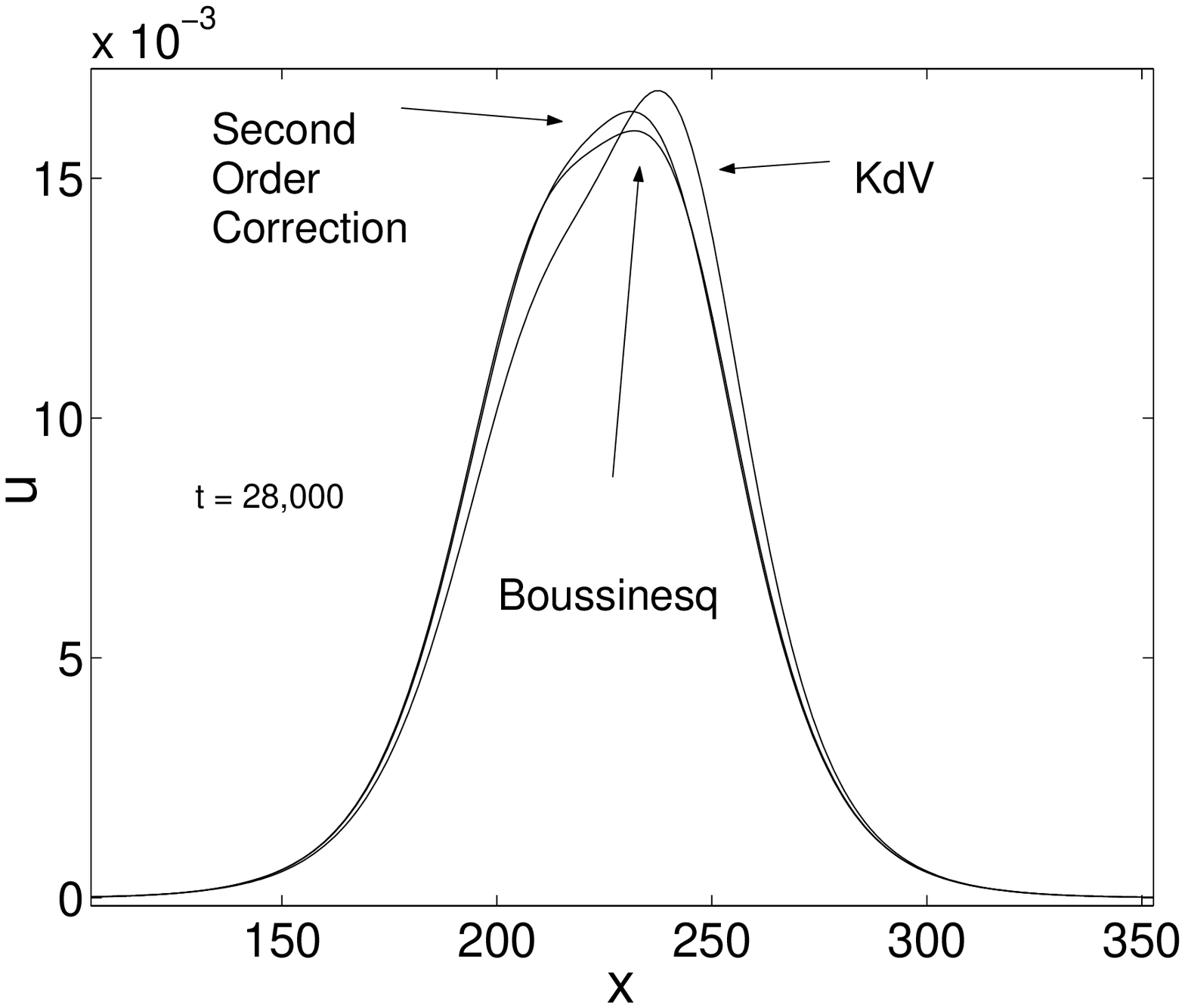, height=2.5in, clip=}
\caption{Overtaking Wave at the Collision, $\epsilon = 0.05$}
\label{overatcollision}
\end{center}
\end{figure}

\begin{figure}
\begin{center}
\epsfig{file=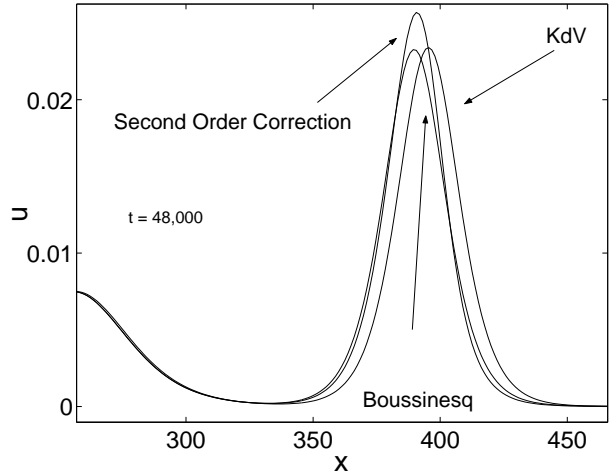, height=2.5in, clip=}
\caption{Overtaking Wave after the Collision, $\epsilon = 0.05$}
\label{overaftercollision}
\end{center}
\end{figure}

\begin{figure}
\begin{center}
\epsfig{file=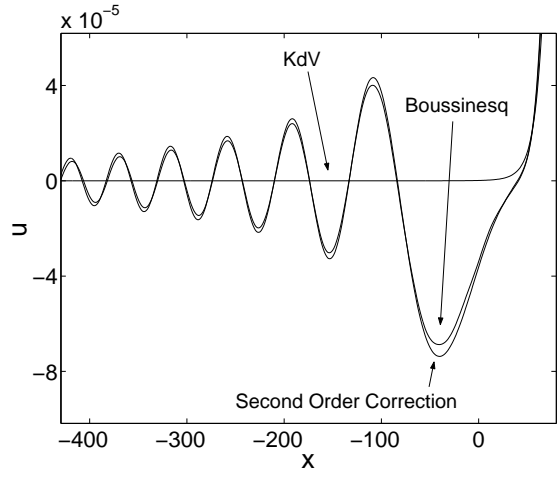, height=2.5in, clip=}
\caption{Dispersive Wave for the Overtaking Wave}
\label{overtakingdwave}  
\end{center}
\end{figure}

Figures \ref{overbeforecollision}, \ref{overatcollision}  and
\ref{overaftercollision} display the values of $u$ and the approximations
at various times during the collision.  
As in the case of the head-on collision, the second order correction
picks up the presence of a dispersive wave, which is not seen in the KdV
approximation (see Figure \ref{overtakingdwave})

\begin{figure}
\begin{center}
\epsfig{file=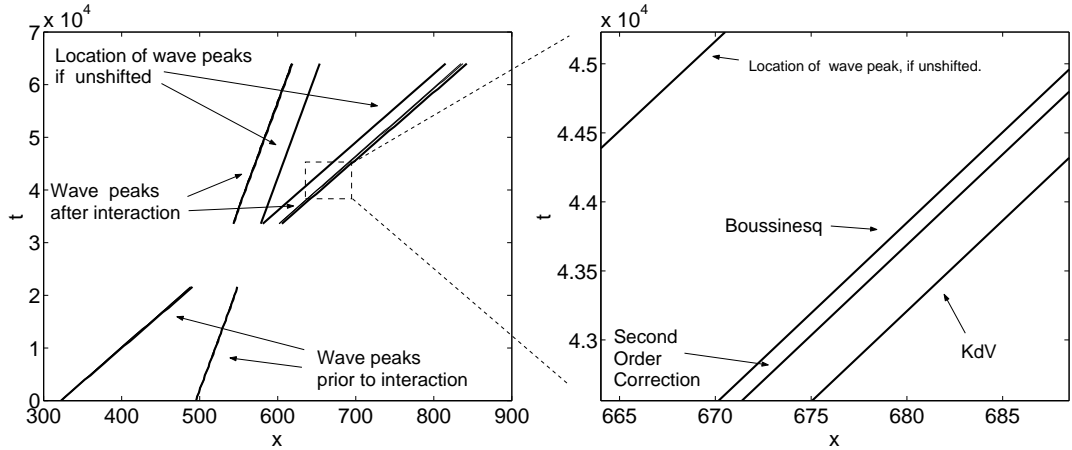, height=2.5in, clip=}
\caption{Wave Peak Locations, $\epsilon = 0.05$}
\label{peaklocations}
\end{center}
\end{figure}

It is well known that in the two soliton interaction, the waves
are phase-shifted after the collision (that is, the faster wave
is further ahead after the collision than it would have been had
no interaction taken place, and the slower wave falls behind in
a similar fashion).  Overtaking waves in the Boussinesq equation
share this feature, though with a different phase shift.  
This can be seen
in Figure \ref{overaftercollision}, where the KdV approximation is
leading the Boussinesq solution.  The second order correction noticeably
``fixes'' this problem.  In Figures \ref{peaklocations}
we plot the locations of the peaks.  Note
that these figures reflect the fact that
the numerics are computed in a moving reference frame (moving
to the right with unit velocity).

\begin{figure}
\begin{center}
\epsfig{file=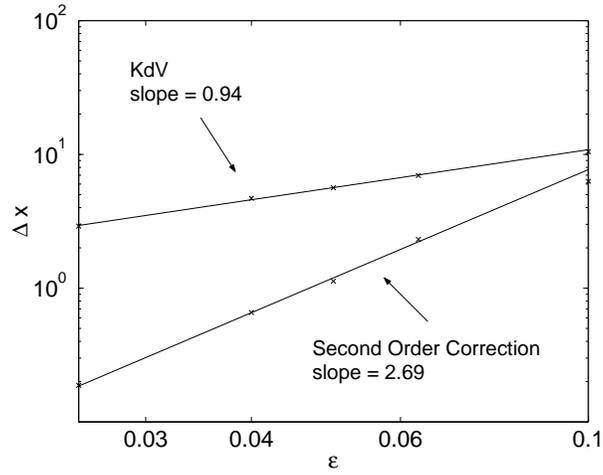,height=2.5in, clip=}
\caption{Error in Phase Shift vs. $\epsilon$}
\label{phaseshiftvsep}
\end{center}
\end{figure}

In Figure \ref{phaseshiftvsep}
we plot the error in the phase shifts for the two approximations
versus $\epsilon$.  Notice that the slope for the second order correction
is steeper than that of the KdV approximation

\begin{figure}
\begin{center}
\epsfig{file=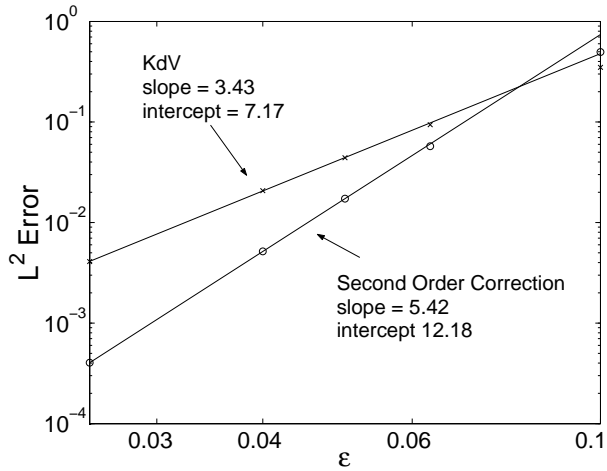, height=2.5in, clip=}
\caption{$\sup \|u - w\|_{L^2}$ vs. $\epsilon$ for Overtaking Wave Collision}
\label{ol2errvsep}
\end{center}
\end{figure}

\begin{figure}
\begin{center}
\epsfig{file=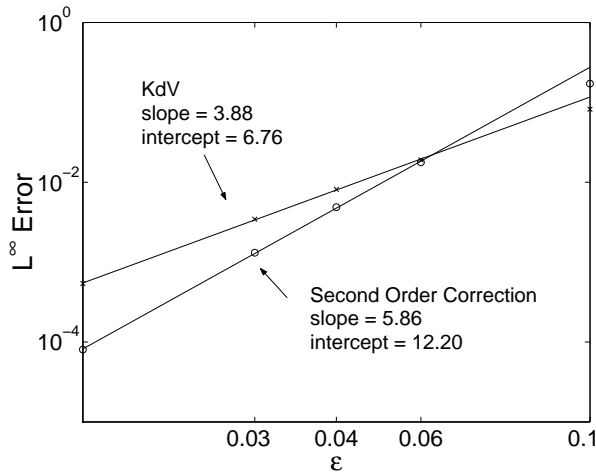, height=2.5in, clip=}
\caption{$\max \|u - w\|_{L^\infty}$ vs. $\epsilon$ for Overtaking Wave Collision}
\label{olinferrvsep}  
\end{center}
\end{figure}

In Figures \ref{ol2errvsep} and \ref{olinferrvsep}
we plot the maximum of the $L^2$ and $L^\infty$ error
for the two approximations versus $\epsilon$ on a log-log
plot (as we did for the head-on interaction earlier).
We summarize the results in Tables (\ref{overtakingslopetable})
and (\ref{overtakingconstanttable}).  

\begin{table}
\begin{tabular}{|l||c|c|} \hline
                              & $L^2$  & $L^\infty$ \\
\hline\
KdV                           & $3.43$ &  $3.88$\\
\hline
KdV + second order correction & $5.42$ &   $5.86$\\
\hline
\end{tabular}
\caption{Order of the approximation, numerically computed, for the overtaking collision}
\label{overtakingslopetable}
\end{table}

\begin{table}
\begin{tabular}{|l||c|c|} \hline
                              & $L^2$  & $L^\infty$ \\
\hline\
KdV                           & $1300$ &  $860$\\
\hline
KdV + second order correction & $19,500$ &   $19,900$\\
\hline
\end{tabular}
\caption{Value of $C_F$, numerically computed, for the overtaking collision}
\label{overtakingconstanttable}
\end{table}

\section{Conclusions}

We conclude by briefly surveying other work on the derivation of
higher order modulation equations for water waves and related
systems.

For the actual water wave equations there have been a number of
studies of corrections to the KdV approximation to water waves
spanning the spectrum from non-rigorous asymptotic expansions
\cite{byatt-smith:71}, \cite{byatt-smith:88},
\cite{su.etal:80} to numerical solutions of the equations
of motion and comparison with the KdV predictions
\cite{bona.etal:81},\cite{su.etal:82}, \cite{fenton.etal:82} to experimental
investigations \cite{maxworthy:76}, \cite{cooker.etal:97}.
We concentrate here on the theoretical studies since they have
the closest connection to our work.  In the investigations
of Byatt-Smith \cite{byatt-smith:71}, \cite{byatt-smith:88} and
Su \& Marie \cite{su.etal:80} the focus is on the head-on collision
of solitary waves.  This has several consequences.  First of
all the authors assume that the initial conditions are of
a special form, namely a pair of counter-propagating
solitary waves.  The higher-order corrections to the solution
then exploit this special form by including not only a
correction to the amplitude of the solution, but a phase
shift for each wave as it undergoes the collision.  This is
a very reasonable hypothesis in these physical circumstances, but one 
which can't easily be adapted to the more general
type of initial conditions considered in our work.
Furthermore, since these papers consider specifically the head-on
collision of solitary waves they are concerned with events
which occur on relatively short time scales (i.e. time
scales of $O(\frac{1}{\epsilon})$ in our scaling.)
As noted in (\cite{byatt-smith:88}, p. 503) these expansions
are not uniformly valid in time and it is not clear
whether or not their solutions could be controlled over
time scales of  $O(\frac{1}{\epsilon^3})$.  It is
worth noting that in spite of the differences
between our approach and those discussed here, Byatt-Smith
\cite{byatt-smith:88} also finds that corrections
to the amplitude of the solitary wave evolve according
to the linearized KdV equation.

Another set of papers by Sachs, \cite{sachs:84}, Zho and Su 
\cite{zho.etal:86}, and H{\u{a}}r{\u{a}}gu{\c{s}}-Courcelle, Sattinger,
and Nicholls \cite{haragus.etal:02} considers corrections
to the KdV approximation for unidirectional motion.  The
first two of these papers study this question in the 
context of water waves, while \cite{haragus.etal:02} studies
the KdV approximation to solutions of the Euler-Poisson
equations.  The focus of these papers (particularly \cite{sachs:84}
and \cite{haragus.etal:02}) is rather different than ours, however.
Both derive an inhomogeneous linearized 
KdV equation for the correction to the KdV approximation.
However, rather than deriving rigorous estimates 
of the difference between
the approximate solutions provided by the model equations
and the true solutions they focus on
the nature of the solutions of the linearized, inhomogeneous KdV equation.
In particular, Sachs \cite{sachs:84}
shows that if one linearizes about the
$N$-soliton solution of the KdV equation the resulting inhomogeneous
equation has solutions which have no secular growth.
In \cite{haragus.etal:02} the authors 
obtain explicit solutions of the linearized KdV equation,
particularly for the case in which one linearizes about
the two-soliton solution of the KdV equation.
H{\u{a}}r{\u{a}}gu{\c{s}}-Courcelle, et al. 
then compare the approximation they obtain to numerically
computed solutions of the Euler-Poisson equation and
they find that the addition of the solution of the linearized
KdV equation to the approximation given by the 
two-soliton solution of the KdV equation does a significantly
better job of approximating the solution of the Euler-Poisson
equation.  In particular, they note that the prediction
of the phase shift that occurs when a ``fast'' traveling
wave overtakes a slower one is significantly better when
the second order correction is included.  This effect is
also present in our 
approximation -- see Figure (\ref{phaseshiftvsep}).

\noindent
{\bf Acknowledgments:}
The work of the authors was supported
in part by the NSF grant DMS-0103915.  The authors thank
D. Nicholls and D. Sattinger for several discussions of
their work \cite{haragus.etal:02}, prior to publication.  The second
author also thanks G. Schneider for numerous discussions
of the derivation and justification of modulation equations and
H. Segur for discussions of higher order approximations to
water waves.

\section{Appendix}
\begin{proof}  {\bf For Lemma \ref{cross-term}}:
We consider the case with the ``minus'' sign on the left hand sign for simplicity.  The other case is
analogous.  First we change variables to $X_+ = X+\tau$, $T=\epsilon^2 \tau$, and $\upsilon(X_+,T)=u(X,\tau)$.  Under this
change we get the equation:
\begin{equation}
\label{new coord trans}
\partial_T \upsilon(X_+,T)=\epsilon^{-2} l(X_+,T)r(X_+ - 2T \epsilon^{-2},T)
\end{equation}
This can solved by integrating with respect to the variable $T$.  We get,
$$
\upsilon(X_+,T)=\epsilon^{-2} \int_0^T l(X_+,s)r(X_+-2s\epsilon^{-2}2,s)ds
$$
Now we multiply by the appropriate weight and take norms:
\begin{align*}
  &(1+X_+^2) |\upsilon(X_+,T)| \\
  \le& \epsilon^{-2} \int_0^T (1+X_+^2) |l(X_+,s)| |r(X_+ - 2s \epsilon^{-2},s)|\ ds\\
  \le& \epsilon^{-2} \int_0^T (1+X_+^2) |l(X_+,s)| (1+(X_+ - 2s \epsilon^{-2})^2) |r(X_+ - 2s \epsilon^{-2},s)|\ ds\\
    =& \epsilon^{-2} \int_0^T \frac{(1+X_+^2)^2 |l(X_+,s)| (1+(X_+ - 2s \epsilon^{-2})^2)^2 |r(X_+ - 2s \epsilon^{-2},s)| }
                                   {(1+X_+^2) (1+(X_+ - 2s \epsilon^{-2})^2)}\ ds\\
    \le& \epsilon^{-2} \int_0^T \frac{(1+X_+^2)^2 |l(X_+,s)| (1+(X_+ - 2s \epsilon^{-2})^2) |r(X_+ - 2s \epsilon^{-2},s)| }
                               {(1 +  (2s \epsilon^{-2})^2)}\ ds
\end{align*}

Now take the $H^s$ norm of each side of this equation and find that:
\begin{align*}
\|\upsilon(\cdot,T)\|_{H^s(2)} 
      &\le \|l\|_{H^s(4)} \|r\|_{H^s(4)} \epsilon^{-2} \int_0^T \frac{1}{(1+(2s \epsilon^{-2})^2)}ds\\
      &\le C\|l\|_{H^s(4)} \|r\|_{H^s(4)} \arctan(2T \epsilon^{-2})\\
      &\le C\|l\|_{H^s(4)} \|r\|_{H^s(4)} 
\end{align*}
\end{proof}

\begin{proof} {\bf For Lemma \ref{inner prod 1}}:
\begin{align*}
\langle u \partial_x f,& \rangle_{H^s(2)} 
 =  \langle u \partial_x f,f \rangle_{H^{s-1}(2)} + ((1+x^2) \partial^s_x(u \partial_x f),(1+x^2) \partial^s_x f)_{L^2} \\
\le& |u \partial_x f|_{H^{s-1}(2)} |f|_{H^{s-1}(2)} 
      + \sum^{s-1}_{j=0} c_{sj}( (1+x^2)^2 \partial^{s-j}_x u \partial^{j+1}_x f,\partial^s_x f)_{L^2}\\
    &+ \int (1+x^2)^2 u \partial_x^{s+1} f \partial^s_x f dx\\
\le& C |u|_{H^s(2)} |f|^2_{H^{s}(2)} + \frac{1}{2}\int (1+x^2)^2 u \partial_x(\partial^s_x f)^2 dx\\
\le& C |u|_{H^s(2)} |f|^2_{H^{s}(2)} - \frac{1}{2}\int \partial_x((1+x^2)^2 u)(\partial^s_x f)^2 dx\\
\le& C |u|_{H^s(2)} |f|^2_{H^{s}(2)} - \frac{1}{2}\int \partial_x u ((1+x^2)\partial^s_x f)^2 dx\\
   & - 2\int x u (1+x^2) (\partial_x^s f)^2 dx\\
\le& C |u|_{H^s(2)} |f|^2_{H^{s}(2)}
\end{align*}
\end{proof}

\begin{proof}  {\bf For Lemma \ref{inner prod 2}}:
In this proof we use the standard norm on $H^s(2)$.
\begin{align*}
(f,\partial^3_x f)_{H^s(2)} &=    -6(x^2 f, \partial_x f)_{H^s} - 6((1+x^2)f,x \partial^2_x f)_{H^s}\\
		 	    &\le  C(\|f\|_{H^s(2)} \|f\|_{H^{s+1}}) + C (\|f\|_{H^s(2)} \|x \partial^2_x f\|_{H^s})
\end{align*}
So now consider 
\begin{equation*}
\|x \partial^2_x f\|^2_{H^s} 
   = \|x \partial^2_x f\|^2_{H^{s-2}} + 
      \| \partial_x^{s-1} (x \partial^2_x f) \|^2_{L^2} + 
      \|\partial^s_x(x \partial^2_x f)\|^2_{L^2}
\end{equation*}
We now treat the last term in the above, as the middle term can be handled in a similar fashion, and the
first is easily dealt with.
\begin{align*}
\|\partial^s_x(x \partial^2_x f)\|^2_{L^2} 
   \le& C \|\partial^{s+1}_x f\|^2_{L^2} + \|x \partial^{s+2}_x f \|_{L^2}\\
   \le& C \|f\|^2_{H^{s+4}} + \int x^2 \partial^{s+2}_x f \partial^{s+2}_x f dx\\
   \le& C \|f\|^2_{H^{s+4}} + \int x^2 \partial^{s+4}_x f \partial^{s}_x f dx 
                          + 4 \int x \partial^{s+3}_x f \partial^{s}_x f dx\\ 
      &+ 2 \int   \partial^{s+2}_x f \partial^{s}_x f dx\\
   \le& C (\|f\|^2_{H^{s+4}} +  \|f\|_{H^{s+4}}\|f\|_{H^{s}(2)})
\end{align*}
This estimate completes the proof.
\end{proof}

\begin{proof} {\bf For Lemma \ref{op-est 2}}:
Notice that the polynomials $T_j$ are the first, third and fifth order polynomial expansions of $y/\sqrt{1-y^2}$
about $y=0$.  Moreover, note that only odd powers appear in the expansion.  So, by Taylor's theorem, 
there is a constant $C$ such that $|y/\sqrt{1-y^2} - T_j(y)| \le C |y|^{j+2}$.

We shall now use the Fourier transform version of the Sobolev norms in the following computation, which 
concludes the proof.  Consider

\begin{eqnarray*}
&&\| \lambda \Phi(\epsilon \cdot)-T_j(\epsilon \partial_X) \Phi(\epsilon \cdot)\|^2_{H^s}
\\&=& 
\int(1+k^2)^s | \widehat{ \lambda \Phi(\epsilon x)}- \widehat{ T_j(\epsilon \partial_X) \Phi(\epsilon x) } |^2 dk\\ 
&=&
\epsilon^{-2} \int(1+k^2)^s | (\frac{ik}{\sqrt{1+k^2}}-T_j(ik))\hat{\Phi}(k/\epsilon) |^2 dk\\ 
&\le&
C \epsilon^{-2} \int(1+k^2)^s | k^{j+2} \hat{\Phi}(k/\epsilon) |^2 dk, \quad K=k/\epsilon\\ 
&=&
C \epsilon^{2j+3} \int(1+(\epsilon K)^2)^s | K^{j+2} \hat{\Phi}(K) |^2 dK\\ 
&\le&
C \epsilon^{2j+3} \int(1+ K^2)^s | K^{j+2} \hat{\Phi}(K) |^2 dK\\ 
&\le&
C \epsilon^{2j+3} \|\partial_X^{j+2} \Phi\|^2_{H^s}\\ 
&\le&
C \epsilon^{2j+3} \|\Phi\|^2_{H^{s+j+2}}
\end{eqnarray*}
\end{proof}

\begin{proof}{\bf For Lemma \ref{following fact}}:
\begin{eqnarray*}
& & \left(R^2-R^1,\lambda[(\Psi^1 +\Psi^2)(R^1+R^2)] \right)_{H^s}\\
&=&\left(-\lambda(R^2-R^1),(\Psi^1 +\Psi^2)(R^1+R^2) \right)_{H^s}\\
&\le&-\left(\partial_t(R^1+R^2),(\Psi^1 +\Psi^2)(R^1+R^2) \right)_{H^s}\\
&+&\epsilon^{-{11/2}}\left(\textrm{Res}[\bar{\Psi}]^1+\textrm{Res}[\bar{\Psi}]^2,(\Psi^1 +\Psi^2)(R^1+R^2)\right)_{H^s}\\
&\le& - \left( \partial_t (R^1+ R^2),(\Psi^1 + \Psi^2) (R^1 + R^2) \right)_{H^s} 
 + C \epsilon^{3} \|\bar{R}\|_{H^s \times H^s}
\end{eqnarray*}
\end{proof}

\begin{proof}{\bf For Lemma \ref{commuter est}}:
We shall be using the fact that, by the Sobolev embedding theorem, we have  $\gamma$ and its first $s$ derivatives in $L^\infty$.  
\begin{align*}
 & | (f(x),\gamma(\epsilon x) g(x))_{H^s}-(g(x),\gamma(\epsilon x) f(x))_{H^s}| \\ 
=& | \sum^s_{j=0} \left(\partial^j_x f(x),\partial^j_x (\gamma(\epsilon x) g(x))\right)_{L^2}
  -\left(\partial^j_x g(x),\partial^j_x(\gamma(\epsilon x) f(x))\right)_{L^2} | \\
=& | \sum^s_{j=0} \sum^j_{l=0} c_{jl} \left\{ \left(\partial^j_x f(x),\partial^l_x (\gamma(\epsilon x)) 
   \partial^{j-l}_x g(x))\right)_{L^2}
  -\left(\partial^j_x g(x),\partial^l_x (\gamma(\epsilon x)) \partial^{j-l}_x f(x))\right)_{L^2} \right\} |\\
=& | \sum^s_{j=0} \sum^j_{l=0} c_{jl} \left(\partial^j_x f(x) \partial^{j-l}_x g(x) 
   - \partial^j g(x) \partial^{j-l}_x f(x),\partial^l_x (\gamma(\epsilon x))\right)_{L^2} |\\
=& \epsilon | \sum^s_{j=0} 
            \sum^j_{l=1} c_{jl} \left(\partial^j_x f(x) \partial^{j-l}_x g(x) - \partial^j g(x) \partial^{j-l}_x f(x)
                                     ,\epsilon^{l-1} \partial^l_X \gamma(\epsilon x)\right)_{L^2} |\\
=& \epsilon | \sum^s_{j=0} 
            \sum^j_{l=1} c_{jl} \int \left(\partial^j_x f(x) \partial^{j-l}_x g(x) - 
            \partial^j g(x) \partial^{j-l}_x f(x)\right)
            \left(\epsilon^{l-1} \partial^l_X \gamma(\epsilon x)\right) dx|\\
\le& \epsilon \|\gamma\|_{W^{s,\infty}} 
             \sum^s_{j=0} 
             \sum^j_{l=1} c_{jl} \int | \left(\partial^j_x f(x) \partial^{j-l}_x g(x)
           - \partial^j g(x) \partial^{j-l}_x f(x)\right)| dx\\
\le& C \epsilon  \|\gamma\|_{W^{s,\infty}} \|f\|_{H^s} \|g\|_{H^s}
\end{align*}
\end{proof}

\begin{proof}{\bf For Lemma \ref{gron}}:
Functions which obey the inequality are bounded above by solutions to the following family
of ordinary differential equations,
$$
\dot{\eta}(T;\epsilon) = C (1 + \eta(T;\epsilon) + \epsilon^{5/2} \eta^{3/2}(T;\epsilon)),\quad \eta(0;\epsilon)=0,
$$
and so we prove the result for these equations.  

By separation of variables, we have that $\eta(T;0)=e^{CT}-1$.  We notice that 
for fixed $T$, $\eta(T;\epsilon)$ is a continuous and increasing function of
$\epsilon$.  This follows since solutions of ODEs depend smoothly on their
parameters and that the right hand side of the differential equations is
increasing in $\epsilon$.  

Thus, by the intermediate value theorem, there
exists $\epsilon_0$ such that $\eta(T_0;\epsilon_0)=\epsilon_0^{-5}$.  Moreover
since $\epsilon^{-5}$ is a decreasing function for $\epsilon>0$, we have
$\eta(T_0;\epsilon)\le\epsilon^{-5}$ for $\epsilon \in (0,\epsilon_0)$.
We further note that for fixed $\epsilon$, $\eta$ is continuous and increasing
in $T$.  So we have $\eta(T;\epsilon) \le \epsilon^{-5}$ 
for $T \in [0,T_0]$
and $\epsilon \in (0,\epsilon_0)$.  

Thus we have,
$$
\dot{\eta}(T;\epsilon) \le C (1 + 2\eta(T;\epsilon))\quad \eta(0;\epsilon)=0,
$$
for $T \in [0,T_0]$
and $\epsilon \in (0,\epsilon_0)$.   We apply Gronwall's inequality
to this to prove the result.
\end{proof}

\nocite{*}

\bibliography{ref2}
\bibliographystyle{plain}

\end{document}